\documentclass[11pt,a4paper,reqno]{amsart}
\usepackage[T1]{fontenc}
\usepackage{lmodern}
\usepackage{amsmath,amssymb,amsthm,mathtools,mathrsfs}
\usepackage[margin=28mm]{geometry}
\usepackage{microtype}
\usepackage{enumitem}
\usepackage{xcolor}
\usepackage{tikz}
\usepackage{tikz-cd}

\usepackage[breaklinks,colorlinks=false,
	citebordercolor={0 1 0},linkbordercolor={1 0 0},
	pdfborder={0 0 1}]{hyperref}
\usepackage{bookmark}
\hypersetup{pdftitle={Odd-degree cycles and rationality of Peskine sixfolds of discriminant 24},pdfauthor={Yusuke Ohmiya}}
\setlist[enumerate]{label=\textup{(\arabic*)},leftmargin=2.1em,itemsep=3pt,topsep=5pt}
\setlist[itemize]{leftmargin=1.8em,itemsep=3pt,topsep=5pt}
\newtheorem{theorem}{Theorem}[section]
\newtheorem{lemma}[theorem]{Lemma}
\newtheorem{proposition}[theorem]{Proposition}
\newtheorem{corollary}[theorem]{Corollary}
\theoremstyle{definition}
\newtheorem{definition}[theorem]{Definition}
\newtheorem{assumption}[theorem]{Assumption}
\newtheorem{conjecture}[theorem]{Conjecture}
\theoremstyle{remark}
\newtheorem{remark}[theorem]{Remark}
\theoremstyle{plain}
\numberwithin{equation}{section}

\newcommand{\rBr}{\operatorname{Br}}
\newcommand{\rk}{\operatorname{rk}}

\newcommand{\Spec}{\operatorname{Spec}}
\newcommand{\Gr}{\operatorname{Gr}}
\newcommand{\Sym}{\operatorname{Sym}}

\newcommand{\RHom}{\operatorname{\mathbf{R}Hom}}
\newcommand{\End}{\operatorname{End}}

\newcommand{\Perf}{\operatorname{Perf}}
\newcommand{\Ku}{\operatorname{Ku}}

\newcommand{\Db}{D^{\mathrm b}}
\newcommand{\Dqc}{D_{\mathrm{qc}}}
\newcommand{\cO}{\mathcal O}
\newcommand{\PP}{\mathbf P}

\newcommand{\CC}{\mathbf C}
\newcommand{\ZZ}{\mathbf Z}
\newcommand{\QQ}{\mathbf Q}
\newcommand{\HH}{\mathbb H}

\title[Odd-degree cycles and rationality of Peskine sixfolds of discriminant 24]{Odd-degree cycles and rationality of Peskine sixfolds of discriminant 24}
\author{Yusuke Ohmiya}

\begin{document}

\begin{abstract}
Benedetti and Faenzi conjectured that, for a general trivector $\sigma\in\mathcal D^{1,6,10}$, the Peskine sixfold $X_1^\sigma$ is rational if the Brauer class $\beta$ on the associated K3 surface vanishes. We prove this implication for every $\sigma$ in a nonempty Zariski open subset $U\subset\mathcal D^{1,6,10}$. When $\beta=0$, we construct a line on the generic fiber of the Peskine fibration after a cubic field extension and use odd-degree descent to obtain a line over the original function field. We extend criteria for rationality to smooth complete intersections of an arbitrary number of quadrics under suitable hypotheses. Finally, we describe the derived category of the generic fiber of the Peskine fibration in terms of a twisted nodal curve of arithmetic genus two.
\end{abstract}
\maketitle

\section{Introduction}
The rationality problem for algebraic varieties is a classical problem that remains central to algebraic geometry. 
Clemens--Griffiths and Artin--Mumford gave distinct obstructions to rationality for unirational threefolds \cite{CG72,AM72}.

A conjecture of Benedetti and Faenzi relates the rationality of Peskine sixfolds on the discriminant $24$ divisor $\mathcal D^{1,6,10}$ to the vanishing of a Brauer class on the associated K3 surface \cite[Conjecture~4.1]{BF24}. The key construction relates a rational section of the Brauer--Severi conic on the associated K3 surface to a line, defined after a cubic field extension, on the generic fiber of the Peskine fibration constructed by Benedetti and Faenzi. Springer's theorem and Amer's theorem then give a line over the original function field. We also formulate relative odd-degree cycle criteria for intersections of quadrics and study the bounded derived categories of coherent sheaves on the related varieties.

\subsection{Background}
\subsubsection{Cubic fourfolds}
Cubic fourfolds provide an important setting for this problem. Hassett studied the divisors $\mathcal C_d$ of special cubic fourfolds, where $d$ is the discriminant of a rank-two lattice of algebraic classes containing the square of the hyperplane class. He determined when an associated polarized K3 surface exists at the level of Hodge structures \cite{Has00} and constructed rational cubic fourfolds containing a plane using odd intersection numbers \cite{Has99}. Kuznetsov proposed a categorical rationality conjecture \cite[Conjecture~1.1]{Kuz10}: a smooth complex cubic fourfold $Y$ is rational if and only if its component $\Ku(Y)$ in
$$
	\Db(Y)=\langle\Ku(Y),\cO_Y,\cO_Y(1),\cO_Y(2)\rangle
$$
is equivalent to $\Db(S)$ for an untwisted projective K3 surface $S$. Addington--Thomas compared the Hodge-theoretic and categorical conditions and proved their agreement on a nonempty Zariski open subset of each relevant divisor \cite{AT14}.

For a smooth cubic fourfold $Y$ containing a plane, projection from the plane gives a quadric surface fibration from the blowup of $Y$ along that plane to $\PP^2_{\CC}$. If its discriminant sextic is smooth, this fibration determines a degree-two K3 surface $S$ and a Brauer class $\alpha\in\rBr(S)[2]$, with $\Ku(Y)\simeq\Db(S,\alpha)$ \cite[Theorem~4.3]{Kuz10}. The condition $\alpha=0$ implies rationality \cite{Has99,ABBV14}. This gives a concrete relation between rationality, an associated K3 surface, and the vanishing of a Brauer class.

\subsubsection{Peskine sixfolds and their rationality}
A similar picture is expected for Peskine sixfolds. Let $V_{10}$ be a ten-dimensional complex vector space. A trivector $\sigma\in\bigwedge^3V_{10}^{\vee}$ defines
$$
	\begin{aligned}
		X_1^\sigma&=\{[U_1]\in\PP(V_{10}):
		\rk\sigma(U_1,V_{10},V_{10})\leq6\},\\
		X_6^\sigma&=\{[U_6]\in\Gr(6,V_{10}):
		\sigma|_{\wedge^3U_6}=0\}.
	\end{aligned}
$$
For general $\sigma\in\bigwedge^3V_{10}^{\vee}$, these are the smooth Peskine sixfold and the Debarre--Voisin hyperk\"ahler fourfold \cite{DV10}. Their relation is analogous to that between a cubic fourfold and its variety of lines \cite{BD85}; see also \cite{BFM21} for related varieties of K3 type. When both $X_1^\sigma$ and $X_6^\sigma$ are smooth, Benedetti--Song identified the vanishing integral cohomology $H^6(X_1^\sigma,\ZZ)_{\mathrm{van}}$, up to a Tate twist, with $H^2(X_6^\sigma,\ZZ)_{\mathrm{prim}}$ \cite[Theorem~1.1]{BS21}. They also studied the divisors $\mathcal D^{3,3,10}$, $\mathcal D^{1,6,10}$, and $\mathcal D^{4,7,7}$ in the moduli space of trivectors $\PP(\bigwedge^3V_{10}^{\vee})/\!/\mathrm{SL}(V_{10})$.

These correspond to Heegner divisors of lattice discriminants $22$, $24$, and $28$, respectively \cite[\S2.2 and Table~1]{BS21}. For general $\sigma\in\mathcal D^{1,6,10}$ they constructed a degree-six K3 surface $S$ and a Brauer--Severi conic $\pi:D\to S$, with class $\beta\in \rBr(S)[2]$ \cite[Propositions~5.5 and 5.7]{BS21}. The connection with cubic fourfolds is further developed in \cite{BFFM26}, which constructs an associated cubic fourfold for a very general trivector on this discriminant $24$ locus. Benedetti and Faenzi proved that $X_1^\sigma$ is rational for general $\sigma\in\mathcal D^{3,3,10}$ \cite[Theorem~2.1]{BF24}. For $\mathcal D^{1,6,10}$, they obtained a rationality criterion and formulated the following conjecture.

\begin{conjecture}[{\cite[Conjecture~4.1]{BF24}}]
	Let $\sigma \in \mathcal{D}^{1,6,10}$ be a general element. If the Brauer class $\beta$ of the associated degree-six twisted K3 surface $(S,\beta)$ vanishes, then the Peskine sixfold $X_1^\sigma$ is rational.
\end{conjecture}

By \cite[Propositions~3.1, 3.2, 3.10, and 3.17--3.19]{BF24}, we can choose a nonempty Zariski open subset $U\subset\mathcal D^{1,6,10}$ as in Assumption~\ref{pes:generic}. For every $\sigma\in U$, Theorem~\ref{intro:peskine} shows that $X_1^\sigma$ is rational if $\beta=0$. Moreover, Proposition~\ref{pes:nonvacuous} shows that the intersection of $U$ with the locus where $\beta=0$ is nonempty. The setup follows \cite[Propositions~3.1, 3.2, 3.10 and 3.17--3.19]{BF24}; we recall it in Section~\ref{pes:geometric-setup}.

For intersections of two quadrics, our relative criterion is a formulation of the quadratic-form argument in Colliot-Th\'el\`ene's appendix to \cite{HT21}. The version for geometrically integral odd-degree curves on intersections with finite singular locus is already established in \cite[Lemma~3.20]{BF24}. We allow cycle classes with integer coefficients without requiring their components on the generic fiber to be geometrically integral, and pass to families by the projection formula. The proof of Theorem~\ref{intro:peskine} uses only the absolute case, namely Theorems~\ref{thm.Springer} and \ref{thm.Amer} and Proposition~\ref{cyc:projection}; the relative criteria place it in the framework of \cite[Theorem~3.21]{BF24}.

For the next theorem, let $k$ be a field of characteristic different from two, let $B$ be a projective $k$-variety of dimension $d$, and let $f:X\to B$ be a dominant projective morphism with $X$ a projective $k$-variety. Let $K=k(B)$ be the function field of $B$. Suppose that $m \geq 2$ and the generic fiber $X_K=X\times_B\operatorname{Spec}  K\subset\PP_K^n$ is an integral complete intersection of $m$ quadrics and that the generic member of their linear system is nondegenerate. Fix an integer $r\geq1$ and a line bundle $\mathcal H$ on $X$ whose restriction to $X_K$ is $\cO_{X_K}(1)$.

\begin{theorem}[Theorem~\ref{cyc:relative}, Corollary~\ref{higher:relative}]\label{intro:relative}
	Assume one of the following conditions:
	\begin{enumerate}
		\item $m=2$, $n\geq4$, and the nonsmooth locus of $X_K\times_K\overline K$ is finite, where $\overline K$ is an algebraic closure of $K$;
		\item $m\geq3$, $n\geq3m-2$, $r\geq m-1$, and $X_K$ is smooth and contains a $K$-linear space $\PP_K^{m-2}$.
	\end{enumerate}
	If there exist a class $z\in\operatorname{CH}_{d+r}(X)$ and line bundles $\mathcal A_1,\ldots,\mathcal A_d$ on $B$ such that
	$$
		\deg_{k}\!\left(
		c_1(\mathcal H)^r
		\prod_{i=1}^{d} f^*c_1(\mathcal A_i)\cap z
		\right)\equiv1\pmod2,
	$$
	then there is a birational map over $B$
	$$
		X\dashrightarrow B\times_{k}\PP_{k}^{n-m}.
	$$
	Consequently, if $B$ is $k$-rational, then $X$ is $k$-rational. In case \textup{(1)}, $X_K$ also contains a $K$-linear space $\PP_K^r$.
\end{theorem}

The usefulness of this formulation is that an odd intersection number on the total space detects a rational generic fiber and, when the base is rational, proves rationality of the total space.

The two cases are proved in Theorem~\ref{cyc:relative} and Corollary~\ref{higher:relative}, respectively. The first follows from the odd-degree argument in \cite[Appendix~A, Theorems~A.3--A.5]{HT21}, in the cycle-theoretic form proved below. The second combines the same argument with standard hyperbolic reduction along the specified linear space; see \cite[Lemmas~3.3 and 3.4]{KS18}. Neither the base nor the total space is required to be smooth.

For the application to Peskine sixfolds, let $\sigma$ belong to the open subset $U\subset\mathcal D^{1,6,10}$ fixed in Assumption~\ref{pes:generic}. Its distinguished flag $V_1\subset V_6\subset V_{10}$ satisfies $\sigma(V_1,V_6,V_{10})=0$. The associated K3 surface $S$, the conic bundle $\pi:D\to S$, and the fibration
$$
	\varpi:X_{1.2}^\sigma\longrightarrow B=\PP(V_{10}/V_6),
$$
where $X_{1.2}^\sigma$ is birational to $X_1^\sigma$, have the properties recalled in Section~\ref{pes:geometric-setup}. The condition $\beta=0$ is not part of this genericity assumption. Write $H$ for the relative hyperplane class and $h=\varpi^*c_1(\cO_B(1))$.

\begin{theorem}[Theorem~\ref{pes:beta-rational}]\label{intro:peskine}
	Let $\sigma\in U$, where $U$ is the open subset fixed in Assumption~\ref{pes:generic}. If $\beta=0$, then the generic fiber of $\varpi$ contains a line $L$ over $\CC(B)$. Its closure $Z_L\subset X_{1.2}^{\sigma}$ satisfies
	$$
		\deg\bigl([Z_L]\cdot H\cdot h^3\bigr)=1.
	$$
	In particular, the Peskine sixfold $X_1^{\sigma}$ is rational.
\end{theorem}

In the proof, we construct the incidence variety $\Gamma$ of degree three, and use it to turn a rational section of $\pi:D\to S$ into compatible projective three-spaces on three members of the Pfaffian pencil.
Intersecting the spaces corresponding to complementary pairs gives a line over a cubic extension of $\CC(B)$. Springer's theorem, applied to this cubic extension, and Amer's theorem then give a line over $\CC(B)$ itself. The closure of this line in $X_{1.2}^\sigma$ satisfies the intersection condition in \cite[Theorem~3.21]{BF24} with intersection number one.

\subsubsection{Derived categories and the Peskine fibration}
The relation between $\Db(X_1^\sigma)$ and $\Db(S,\beta)$ remains to be understood, but we describe a semiorthogonal decomposition of $\Db(X_{1.2/K}^\sigma)$, where $X_{1.2/K}^\sigma$ is the generic fiber of $\varpi:X_{1.2}^\sigma\to B$ and $K=\CC(B)$. In Section~\ref{der:section}, we prove the following result using the decomposition corresponding to the Clifford algebra in \cite{Kuz08Q} and \cite[Theorem~7.25]{CO26}.

\begin{theorem}[Corollary~\ref{der:peskine}]\label{intro:category}
	In the setup of Section~\ref{pes:geometric-setup}, put $K=\CC(B)$ and let $X_{1.2/K}^\sigma$ be the generic fiber of $\varpi$. There exist a projective nodal $K$-curve $M_K$, a class $\gamma_K\in \operatorname{Br}(M_K)$, and a fully faithful functor
	$$
		\Psi_K:\Db(M_K,\gamma_K)\longrightarrow\Db(X_{1.2/K}^\sigma)
	$$
	such that
	$$
		\Db(X_{1.2/K}^\sigma)=\left\langle
		\Psi_K\Db(M_K,\gamma_K),\cO_{X_{1.2/K}^\sigma},\cO_{X_{1.2/K}^\sigma}(1)
		\right\rangle.
	$$
	After base change to an algebraic closure of $K$, the curve $M_K$ has five rational components and six nodes. Its dual graph is the complete bipartite graph $K_{2,3}$, and its arithmetic genus is two.
\end{theorem}

\subsection{Organization of the paper}
Section~\ref{pre:section} collects preliminaries on quadratic forms, relative Fano schemes, Brauer classes, and Peskine varieties. Section~\ref{cyc:section} proves the odd-degree criteria for intersections of two quadrics, and Section~\ref{higher:section} treats complete intersections of an arbitrary number of quadrics by hyperbolic reduction. Section~\ref{cub:section} constructs odd-degree one-cycles from linear spaces on three members of a pencil. Section~\ref{pes:section} proves Theorem~\ref{intro:peskine}, and Section~\ref{der:section} proves Theorem~\ref{intro:category}. Section~8 discusses open questions. Appendix~\ref{app:hyperbolic} gives the coordinate construction and birational maps used in Section~\ref{higher:section}. Appendix~\ref{app:nonvacuous} proves the two lemmas used in Proposition~\ref{pes:nonvacuous}.

\subsection{Notation and conventions}
A variety is an integral separated scheme of finite type over a field. Projective bundles parametrize one-dimensional subspaces: for a vector bundle $\mathcal E$ on a scheme $S$, the bundle $\PP_S(\mathcal E)$ carries the tautological line subbundle $\cO_{\PP_S(\mathcal E)}(-1)\subset p^*\mathcal E$, where $p:\PP_S(\mathcal E)\to S$ is the projection. Let $S$ be a scheme, let $\mathcal E$ be a vector bundle of rank $n$ on $S$, and let $0\leq d\leq n$. The relative Grassmannian $\operatorname{Gr}_S(d,\mathcal E)\longrightarrow S$ is the $S$-scheme characterized by the following property: for every morphism $f:T\to S$, there is a natural bijection
$$
	\operatorname{Hom}_S\bigl(T,\operatorname{Gr}_S(d,\mathcal E)\bigr)
	\cong
	\left\{
	\begin{array}{@{}l@{}}
		\text{locally free }\mathcal O_T\text{-submodules }
		\mathcal F\subset f^*\mathcal E\text{ of rank }d\\
		\text{such that }f^*\mathcal E/\mathcal F
		\text{ is locally free of rank }n-d
	\end{array}
	\right\}.
$$
These bijections are compatible with pullback along morphisms of $S$-schemes. With our convention for projective bundles, we have $\operatorname{Gr}_S(1,\mathcal E)\cong\mathbf P_S(\mathcal E)$.

The notation $\operatorname{CH}_r(Y)$ denotes the Chow group of cycles of dimension $r$ with integer coefficients. For a projective scheme $Y$ over a field $F$, a line bundle $\mathcal H$ on $Y$, and $z\in\operatorname{CH}_r(Y)$, set
$$
	\deg_{\mathcal H}(z)
	=\deg_F\bigl(c_1(\mathcal H)^r\cap z\bigr).
$$
The cycles need not be effective, and their integral components need not be geometrically integral. For a scheme $S$, we write $\Db(S)$ for the bounded derived category of coherent sheaves, $\Dqc(S)$ for the derived category of quasi-coherent sheaves, and $\operatorname{Perf}(S)$ for its subcategory of perfect complexes. For a sheaf of $\mathcal O_S$-algebras $\mathcal A$, we write $D^b(S,\mathcal A)$ for the bounded derived category of coherent right $\mathcal A$-modules, $\Dqc(S,\mathcal A)$ for the derived category of quasi-coherent right $\mathcal A$-modules, and $\operatorname{Perf}(S,\mathcal A)$ for its subcategory of perfect complexes.

\subsection{AI disclosure}
In preparing parts of this manuscript, the author used ChatGPT 6 Astra and Claude Fable 5.1, personally verified all ideas and mathematical content suggested by these models, and reworked the arguments. The author takes full responsibility for the mathematical content of this manuscript.

\section{Preliminaries}\label{pre:section}
\subsection{Some basic notions of quadratic forms and Pfaffians}\label{pre:quadratic-forms}
Let $F$ be a field of characteristic different from two and let $E$ be a finite-dimensional $F$-vector space. For a quadratic form $q:E\to F$, its \textit{polar form} $b_q:E\times E\to F$ is defined by
$$
	b_q(v,w)=q(v+w)-q(v)-q(w).
$$
In particular, $b_q(v,v)=2q(v)$. The {\it rank} and {\it corank} of $q$ are $\operatorname{rk} b_q$ and $\operatorname{corank} b_q$, respectively. A subspace $A\subset E$ is \emph{totally isotropic} if $q|_A=0$; equivalently, $\PP_F(A)$ is contained in the {\it quadric} $V(q)\subset\PP_F(E)$. For nondegenerate $q$, its \emph{Witt index} $i_W(q)$ is the maximum dimension of a totally isotropic $F$-subspace. The {\it radical} of $q$ is $\operatorname{rad}(q)=\{v\in E\mid b_q(v,w)=0\text{ for all }w\in E\}$. Since $b_q(v,v)=2q(v)$, every vector in $\operatorname{rad}(q)$ is isotropic.

Recall that the {\it hyperbolic plane}, denoted by $\mathbb H$, is the two-dimensional $F$-vector space $V=F \oplus F$ equipped with the quadratic form
$$
	\mathbb H : V \longrightarrow F,\qquad (x,y)\longmapsto xy.
$$
The \textit{orthogonal sum} of two quadratic forms $q_1 : V_1 \to F$ and $q_2 : V_2 \to F$ is the quadratic form $q_1 \perp q_2 : V_1\oplus V_2 \to F$, $(v_1,v_2)\mapsto q_1(v_1)+q_2(v_2)$. If $q$ is nondegenerate, Witt decomposition gives
$$
	q\simeq \underbrace{\mathbb H \perp \cdots \perp \mathbb H}_{i_W(q) \textrm{-times}} \perp q_{\mathrm{an}},
$$
where $q_{\mathrm{an}}$ is anisotropic. In particular, splitting off a hyperbolic plane decreases the Witt index by one. Namely, if $q \simeq \mathbb H \perp q'$, then $i_W(q') = i_W(q)-1$. These conventions agree with those of \cite[Appendix~A]{HT21}. We recall the following two results from Colliot-Th\'el\`ene's appendix to \cite{HT21}.

\begin{theorem}[{\cite[Lemma~A.1]{HT21}}] \label{thm.Springer}
	Let $K$ be a field of characteristic different from $2$, and let $L/K$ be a finite extension of odd degree. Let $q:V\to K$ be a quadratic form on a finite-dimensional $K$-vector space $V$. If $q$ is anisotropic over $K$, then its scalar extension
	$$
		q_L:V\otimes_K L\longrightarrow L
	$$
	is anisotropic over $L$.
\end{theorem}

For a field $K$, $K(t)$ denotes the fraction field of the polynomial ring $K[t]$.
\begin{theorem}[{\cite[Theorem~A.4]{HT21}}] \label{thm.Amer}
	Let $K$ be a field of characteristic different from $2$, let $V$ be a finite-dimensional $K$-vector space, and let $q_0,q_1\colon V\to K$ be quadratic forms. Let $t$ be an indeterminate over $K$. After extending scalars to $K(t)$, consider the quadratic form
	$$
		q_0+tq_1\colon V\otimes_K K(t)\longrightarrow K(t).
	$$
	Fix an integer $d$ with $1\leq d\leq \dim_K V$. Then the following conditions are equivalent:
	\begin{enumerate}
		\item There exists a $K(t)$-linear subspace $W\subset V\otimes_K K(t)$ such that
		$$
			\dim_{K(t)}W=d, \qquad (q_0+tq_1)|_W=0.
		$$
		\item There exists a $K$-linear subspace $U\subset V$ such that
		$$
			\dim_K U=d, \qquad q_0|_U=q_1|_U=0.
		$$
	\end{enumerate}
\end{theorem}

If $F$ has characteristic zero, $V$ has dimension $2m$, and $\omega\in\bigwedge^2V^\vee$, a choice of volume form defines the Pfaffian by
$$
	\frac{\omega^{\wedge m}}{m!}
	=\operatorname{Pf}(\omega)\,\mathrm{vol}_V.
$$
For a skew-symmetric matrix $A$ representing $\omega$ in a compatible basis, $\operatorname{Pf}(A)^2=\det(A)$. Thus $\operatorname{Pf}(\omega)=0$ if and only if $\rk\omega\leq2m-2$.

These notions extend to vector bundles as follows. Recall that an integer $n$ is invertible on a scheme $S$ if its image in $\mathcal O_{S,s}$ is a unit for every point $s\in S$. Let $S$ be a scheme on which $2$ is invertible, let $\mathcal{E}$ be a vector bundle on $S$, and let $\mathcal{L}$ be a line bundle on $S$. A {\it quadratic form} on $\mathcal E$ with values in $\mathcal L$ is a morphism of sheaves of sets $q : \mathcal{E} \to \mathcal{L}$ with $q(fv)=f^2q(v)$ for local sections $f$ of $\mathcal O_S$ and $v$ of $\mathcal E$, whose polar form $b_q(v,v') = q(v+v')-q(v)-q(v')$ is $\mathcal{O}_S$-bilinear. Since $2$ is invertible on $S$, such forms correspond bijectively to morphisms $\tilde q : \operatorname{Sym}^2 \mathcal{E} \to \mathcal{L}$.

\subsection{Relative Fano schemes}\label{pre:fano}
We next recall relative parameter spaces of linear subspaces. They can be defined in greater generality, but we only give the special case needed here (for details, see \cite{AK77}, \cite{Kuz14} or \cite{Xie23}).

Let $S$ be a noetherian scheme, let $\mathcal E$ be a vector bundle on $S$, and let $Y\subset\PP_S(\mathcal E)$ be a closed subscheme. For $r\geq0$, let $p : \operatorname{Gr}_S(r+1,\mathcal{E}) \to S$ be the relative Grassmannian. The \emph{relative Fano scheme of projective $r$-planes} is the closed subscheme
$$
	\mathcal{M}_r(Y/S)\subset\operatorname{Gr}_S(r+1,\mathcal E)
$$
whose points over an $S$-scheme $T$ are rank-$(r+1)$ subbundles $\mathcal A\subset\mathcal E_T$ with locally free quotient such that $\PP_T(\mathcal A)\subset Y_T$ as closed subschemes. In particular, $\mathcal{M}_1(Y/S)$ is the relative Fano scheme of lines, whereas $\mathcal{M}_3(Y/S)$ parametrizes projective three-spaces. When $S=\operatorname{Spec} K$ for a field $K$, we also write $\mathcal{M}_r(Y/K)$.

Let $q:\mathcal E\to\mathcal L$ be a quadratic form on $S$, and let $\mathcal Q=V(q)\subset\mathbf P_S(\mathcal E)$ be the quadric it defines. Fix an integer $0\le r<\operatorname{rank}\mathcal E$, and let $p:\operatorname{Gr}_S(r+1,\mathcal E)\longrightarrow S$ be the projection. Then the restriction of $p^*q$ to the tautological subbundle $\mathcal U\subset p^*\mathcal E$ of rank $r+1$ defines a section
$$
	s_q\in H^0\!\left(
	\operatorname{Gr}_S(r+1,\mathcal E),
	\operatorname{Sym}^2\mathcal U^\vee\otimes p^*\mathcal L
	\right),
$$
and hence the relative Fano scheme of $\mathcal Q\to S$ is $\mathcal M_r(\mathcal Q/S)=V(s_q)$. More generally, let $q_i:\mathcal E\to\mathcal L_i$, for $1\le i\le n$, be quadratic forms, and let
$$
	X=V(q_1,\ldots,q_n)\subset\mathbf P_S(\mathcal E)
$$
be their zero scheme. Then
$$
	\mathcal M_r(X/S)=V(s_{q_1},\ldots,s_{q_n})
	\subset\operatorname{Gr}_S(r+1,\mathcal E).
$$
This construction commutes with arbitrary base change on $S$. No smoothness assumption on the quadrics $V(q_i)$ is required, and the resulting relative Fano scheme may be nonreduced.

\subsection{Brauer classes and twisted derived categories}
For a quasi-projective scheme $S$ over a field, we write $\operatorname{Br}(S)=H^2_{\mathrm{\acute et}}(S,\mathbf G_m)_{\rm tors}$ for its Brauer group and $\operatorname{Br}(S)[2] :=\{\alpha\in\operatorname{Br}(S)\mid 2\alpha=0\}$ for its $2$-torsion subgroup. In this setting, the cohomological Brauer group agrees with the group of Morita-equivalence classes of Azumaya $\mathcal O_S$-algebras by a theorem of Gabber; see \cite{de03}, \cite[Tag~0A2J]{Stacks}. For $\alpha\in\operatorname{Br}(S)$, choose an Azumaya algebra $\mathcal A$ representing $\alpha$ and write $D^b(S,\alpha):=D^b(S,\mathcal A)$. This category is independent, up to equivalence, of the choice of $\mathcal A$. We call $D^b(S,\alpha)$ the {\it $\alpha$-twisted derived category} of $S$.

A {\it Brauer--Severi conic} over $S$ is an $S$-scheme $\pi:D\to S$ for which there exists an \'etale covering $\{U_i\to S\}$ and isomorphisms
$$
	D\times_S U_i\simeq\mathbb P^1_{U_i}
$$
over $U_i$. The transition functions of these local trivializations take values in $\mathrm{PGL}_2$. Through the conjugation action on $2\times2$ matrices, they define an Azumaya $\mathcal O_S$-algebra $\mathcal A$ of rank $4$. We call $[\mathcal A]\in\operatorname{Br}(S)$ the Brauer class of $\pi$. Since $\mathcal A$ has degree $2$, its Brauer class satisfies $2[\mathcal A]=0$, and hence belongs to $\operatorname{Br}(S)[2]$.

\subsection{The Peskine varieties}\label{pes:Preliminaries}
We recall Peskine varieties and the relevant moduli spaces; see \cite[\S2.2]{BS21} for details. Let $V_{10}$ be a ten-dimensional complex vector space. A trivector $\sigma\in\bigwedge^3V_{10}^{\vee}$ defines
\begin{equation}\label{def:Peskine}
	X_1^\sigma=\{[U_1]\in\PP(V_{10}) \mid \operatorname{rk}\sigma(U_1,V_{10},V_{10})\leq6\},
\end{equation}
\begin{equation}\label{def:Debarre-Voisin}
	X_6^\sigma=\{[U_6]\in\operatorname{Gr}(6,V_{10}) \mid \sigma|_{\wedge^3U_6}=0\}.
\end{equation}
We call $X_1^\sigma$ a {\it Peskine variety} and $X_6^\sigma$ a {\it Debarre--Voisin variety}.

The GIT quotient
$$
	\mathcal M
	=
	\mathbf P\!\left(\bigwedge^3 V_{10}^{\vee}\right)
	/\!/\operatorname{SL}(V_{10})
$$
is an irreducible projective variety of dimension $20$. Let $\mathcal M_{22}^{(2)}$ denote the moduli space of primitively polarized hyperk\"ahler fourfolds $(X,H)$ of $\mathrm{K3}^{[2]}$-type satisfying
$$
	q(H,H)=22,
	\qquad
	q\bigl(H,H^2(X,\mathbf Z)\bigr)=2\mathbf Z,
$$
where $q$ is the Beauville--Bogomolov--Fujiki form. The construction of \cite{DV10} and the period map give
$$
	\mathfrak m:\mathcal M\dashrightarrow\mathcal M_{22}^{(2)}, \qquad \mathfrak p:\mathcal M_{22}^{(2)}\hookrightarrow\mathcal P.
$$
Here $\mathfrak m([\sigma])=[(X_6^\sigma,H)]$ whenever $X_6^\sigma$ is a smooth fourfold, $H$ denotes its Pl\"ucker polarization, and $\mathcal P$ is the corresponding period space. The map $\mathfrak m$ is birational, and $\mathfrak p$ is an open immersion \cite[Theorem~2.2]{BS21}.

Let $\Lambda$ be the second cohomology lattice of a $\mathrm{K3}^{[2]}$-type fourfold, and fix a primitive element $h\in\Lambda$ of square $22$ and divisibility $2$. For a saturated sublattice $K\subset\Lambda$ of rank $2$ and signature $(1,1)$ containing $h$, the image in $\mathcal P$ of the periods orthogonal to $K$ is a Heegner divisor $\mathcal D_{2e}$, where
$$
	2e=|\operatorname{disc}(K^\perp)|
	=\frac12|\operatorname{disc}(K)|,
$$
and the orthogonal complement is taken in $\Lambda$ \cite[\S2.2 and Lemma~2.1]{BS21}. When nonempty, the inverse image $\mathfrak p^{-1}(\mathcal D_{2e})$ is the corresponding Noether--Lefschetz divisor.

$$
	\begin{tikzcd}[column sep=large, row sep=large]
		& \mathfrak p^{-1}(\mathcal D_{2e})\arrow[r, hook]\arrow[d, hook] & \mathcal D_{2e} \arrow[d, hook] \\
		\mathcal M \arrow[r, dashed, "\mathfrak m"'] & \mathcal M_{22}^{(2)} \arrow[r, hook, "\mathfrak p"'] & \mathcal P
	\end{tikzcd}
$$

The Heegner divisor relevant to this paper is $\mathcal D_{24}$. It corresponds to the divisor $\mathcal D^{1,6,10}\subset\mathcal M$ parametrizing trivectors $\sigma$ that admit a flag $V_1\subset V_6\subset V_{10}$ satisfying
$$
	\dim V_1=1,\qquad \dim V_6=6,\qquad \sigma(V_1,V_6,V_{10})=0.
$$
By \cite[Proposition~5.2]{BS21}, the moduli map $\mathfrak{m}$ restricts to a birational map
$$
	\mathfrak m|_{\mathcal D^{1,6,10}}: \mathcal D^{1,6,10} \dashrightarrow \mathfrak p^{-1}(\mathcal D_{24}).
$$
Thus, for general $\sigma\in\mathcal D^{1,6,10}$, the period of $X_6^\sigma$ lies in $\mathcal D_{24}$. This birationality is used in Proposition~\ref{pes:nonvacuous}.

\section{Odd-degree cycles and a relative rationality criterion}
\label{cyc:section}

Throughout this section, all fields have characteristic different from $2$. We formulate the odd-degree argument of \cite[Appendix~A]{HT21} for cycles with integer coefficients and then apply it to families. For background on intersections of two quadrics over nonclosed fields, see also \cite{CTS87}.

The following is a cycle-theoretic version of \cite[Theorem~A.3]{HT21}. We include a direct proof that does not assume geometric integrality of the components of the cycle.

\begin{lemma}\label{cyc:quadric}
	Let $E$ be a finite-dimensional vector space over $F$, let $q:E\to F$ be a nondegenerate quadratic form, and put $Q=V(q)\subset\PP_F(E)$. Fix $r\geq0$. If $Q$ has an $r$-cycle of odd degree, then the Witt index of $q$ is at least $r+1$. Equivalently, there is an $(r+1)$-dimensional $F$-subspace $A_{r+1}\subset E$ such that $\PP_F(A_{r+1})\subset Q$.
\end{lemma}

\begin{proof}
	Let $z$ be an $r$-cycle of odd degree on $Q$. Suppose first that $F$ is infinite. We argue by induction on $r$. For $r = 0$, write $z = \sum n_i [x_i]$, where the $x_i$ are closed points of $Q$. Since $\deg_F z = \sum n_i[\kappa(x_i):F]$ is odd, the extension $\kappa(x_i) /F$ has odd degree for some $i$. The form $q_{\kappa(x_i)}$ is isotropic, so Theorem~\ref{thm.Springer} implies that $q$ is isotropic over $F$. Now fix $r\geq1$ and assume the assertion for cycles of dimension $r-1$. Since $z$ has odd degree, $z$ has an integral component $W\subset Q$ of dimension $r$ and odd degree. The zero-cycle class $c_1(\mathcal O_Q(1))^r\cap[W]$ has degree $\deg W$, which is odd. The case $r=0$ therefore gives $Q(F)\ne\varnothing$. Since $q$ is nondegenerate, $Q$ is smooth over $F$. The inequality $\operatorname{dim}W=r\ge1$ shows that $Q$ is positive-dimensional, hence geometrically integral of degree $2$. Since $W$ has odd degree, we have $W\ne Q$, so $r<\dim Q$ and $\dim Q\ge2$. Projection from an $F$-point gives a birational map $Q\dashrightarrow\PP_F^{\dim Q}$, and $Q(F)$ is therefore Zariski dense because $F$ is infinite. Write $b_q:E\times E\to F$ for the polar form of $q$. For each $p=[v_p]\in Q(F)$, put $T_p=\{[v]\in\mathbf P_F(E)\mid b_q(v_p,v)=0\}$. We may choose $p\in Q(F)$ such that
	$$
		p\notin W,\qquad W\not\subset T_p.
	$$
	Let $V_W\subset E$ be the smallest vector subspace with $W\subset\mathbb P_F(V_W)$, and let $V_W^\perp$ be its orthogonal complement for $b_q$. The condition $W\subset T_p$ is equivalent to $p\in\mathbb P_F(V_W^\perp)$. Thus the required points lie in
	$$
		Q\setminus\bigl(W\cup(Q\cap\mathbb P_F(V_W^\perp))\bigr).
	$$
	Both $W$ and $Q\cap\mathbb P_F(V_W^\perp)$ are proper closed subsets of $Q$: we have $W\ne Q$, and $Q$ is not contained in any hyperplane. Since $Q$ is integral, the displayed open subset is nonempty. The density of $Q(F)$ gives the required point $p$. Fix a nonzero representative $v_p\in E$ of $p$. The quadratic form
	$$
		\overline q:v_p^\perp/Fv_p\longrightarrow F,
		\qquad v+Fv_p\longmapsto q(v)
	$$
	is well defined and nondegenerate, and $q\simeq\mathbb H\perp\overline q$, where $\mathbb H$ is a hyperbolic plane. Put $\overline Q=V(\overline q)\subset\PP_F(v_p^\perp/Fv_p)$. The hyperplane section $W\cap T_p$ has pure dimension $r-1$ and does not contain $p$. Projection from $p$ therefore defines a proper morphism
	$$
		\pi:W\cap T_p\longrightarrow\overline Q,
		\qquad
		\pi^*\mathcal O_{\overline Q}(1)
		\simeq\mathcal O_{W\cap T_p}(1).
	$$
	Here $[W\cap T_p]$ denotes the fundamental cycle of the hyperplane section, with its intersection multiplicities. The projection formula gives
	$$
		\deg_F\!\left(
		c_1(\mathcal O_{\overline Q}(1))^{r-1}
		\cap\pi_*[W\cap T_p]\right)=\deg W.
	$$
	Thus $\overline Q$ has an $(r-1)$-cycle $\pi_*[W\cap T_p]$ of odd degree with integer coefficients. The induction hypothesis gives $i_W(\overline q)\geq r$. Since $q\simeq\mathbb H\perp\overline q$, we obtain
	$$
		i_W(q)=1+i_W(\overline q)\geq r+1.
	$$
    If $F$ is finite, extend scalars to the infinite field $F(t)$, the fraction field of $F[t]$. The scalar extension $z_{F(t)}$ has the same degree as $z$, so the argument above gives an $F(t)$-point of the projective $F$-scheme
	$$
		\operatorname{OGr}_F(r+1,q)
		\subset \operatorname{Gr}_F(r+1,E)
	$$
	parametrizing totally isotropic subspaces of dimension $r+1$ in $E$. By the valuative criterion of properness, this point extends to a morphism
	$$
		\rho:\operatorname{Spec} F[t]_{(t)}\longrightarrow\operatorname{OGr}_F(r+1,q).
	$$
	Restricting $\rho$ to $\operatorname{Spec}F$ gives an $F$-point of $\operatorname{OGr}_F(r+1,q)$, hence a totally isotropic subspace of dimension $r+1$ in $E$.
\end{proof}

\begin{proposition}\label{cyc:linear-space}
	Let $E$ be a finite-dimensional $K$-vector space, and let $q_0,q_1:E\to K$ be quadratic forms. Suppose that $q_0+tq_1$ is nondegenerate over $K(t)$, and set
	$$
		X=V(q_0,q_1)\subset\PP_K(E).
	$$
	For each fixed $r\geq0$, the following conditions are equivalent:
	\begin{enumerate}
		\item $X$ has an $r$-cycle of odd degree with integer coefficients;
		\item there is an $(r+1)$-dimensional $K$-subspace $A_{r+1}\subset E$ with $\PP_K(A_{r+1})\subset X$.
	\end{enumerate}
\end{proposition}

\begin{proof}
	Assume $(2)$. The cycle $[\PP_K(A_{r+1})]$ has degree one, so $(1)$ follows. Conversely, let $z\in\operatorname{CH}_r(X)$ have odd degree, and let $z_{K(t)}$ be its scalar extension. Put $Q_t=V(q_0+tq_1)$, and let
	$$
		j:X_{K(t)}\hookrightarrow Q_t
		\subset\PP_{K(t)}(E\otimes_K K(t))
	$$
	be the closed immersion. Since $j^*\mathcal O_{Q_t}(1)=\mathcal O_{X_{K(t)}}(1)$, the projection formula gives
	$$
		\deg_{K(t)}\bigl(c_1(\mathcal O_{Q_t}(1))^r\cap j_*z_{K(t)}\bigr)
		=\deg_K\bigl(c_1(\mathcal O_X(1))^r\cap z\bigr).
	$$
	Thus $j_*z_{K(t)}$ has odd degree. Lemma~\ref{cyc:quadric} gives a totally isotropic subspace of dimension $r+1$ for $q_0+tq_1$ over $K(t)$. Applying Theorem~\ref{thm.Amer} with $d=r+1$, we obtain an $(r+1)$-dimensional $K$-subspace $A_{r+1}\subset E$ satisfying
	$$
		q_0|_{A_{r+1}}=q_1|_{A_{r+1}}=0.
	$$
	Hence $\PP_K(A_{r+1})\subset V(q_0,q_1)=X$, which proves $(2)$.
\end{proof}

\begin{remark}
	For smooth intersections and geometrically integral odd-degree subvarieties, the implication from $(1)$ to $(2)$ is \cite[Theorem~A.5]{HT21}. Its proof uses nondegeneracy of the generic quadratic form and Amer's theorem, which also apply to singular base loci. Allowing integer coefficients is formal: an odd-degree cycle has an integral component of odd degree. Lemma~\ref{cyc:quadric} applies without requiring this component to be geometrically integral.
\end{remark}

We now use the preceding results on quadratic forms to obtain a rationality criterion for intersections of two quadrics. For smooth $X$, Proposition~\ref{cyc:projection} is the classical projection from a line \cite[Proposition~2.2]{CTS87}; see also \cite[\S3]{HT21}. We include a proof because $X$ may be singular here.

\begin{proposition}\label{cyc:projection}
	Let $\dim_{K}E=n+1$, where $n\geq4$, and let $X=V(q_0,q_1)\subset\PP_K(E)$ be an integral complete intersection. Assume that $q_0+tq_1$ is nondegenerate over $K(t)$. If $X$ contains a $K$-line $L=\PP_K(A_2)$ that is not contained in the nonsmooth locus of $X\to\operatorname{Spec} K$, then projection from $L$ is a birational map
	$$
		\operatorname{pr}_L:X\dashrightarrow\PP_K(E/A_2)
		\simeq\PP_K^{n-2}.
	$$
\end{proposition}

\begin{proof}
	Since $L=\PP_K(A_2)\subset X$, we have $q_0|_{A_2}=q_1|_{A_2}=0$. Choose coordinates $(x_0,x_1,y := (y_1,\ldots,y_{n-1}))$ on $E$ with $A_2=\{y=0\}$. For $i\in\{0,1\}$, write
	$$
		q_i=x_0\ell_{i0}(y)+x_1\ell_{i1}(y)+Q_i(y),
	$$
	where $\ell_{ij}:E/A_2\to K$ is linear for $j\in\{0,1\}$ and $Q_i:E/A_2\to K$ is quadratic. Put
	$$
		\begin{aligned}
			M(y)&=\begin{pmatrix}\ell_{00}(y)&\ell_{01}(y)\\\ell_{10}(y)&\ell_{11}(y)\end{pmatrix},\\
			\delta(y)&=\det M(y)=\ell_{00}(y)\ell_{11}(y)-\ell_{01}(y)\ell_{10}(y).
		\end{aligned}
	$$
	We first show that $\delta\ne0$. Suppose that $\delta = 0$. If $\ell_{00},\ell_{01}$ are linearly independent, unique factorization gives $\ell_{10} = c \ell_{00}$ and $\ell_{11}=c\ell_{01}$ for some $c\in K$. The quadratic form $q_1-cq_0$ depends only on $y$, so $d(q_1-cq_0)$ vanishes along $L$. The Jacobian matrix of $(q_0,q_1)$ therefore has rank at most one along $L$. This contradicts the assumption that $L$ is not contained in the nonsmooth locus of $X$. If $\ell_{00},\ell_{01}$ are linearly dependent and not both zero, a change of $x_0,x_1$ gives $\ell_{01}=0$ and $\ell_{00}\ne0$. The equality $\delta=0$ forces $\ell_{11}=0$, so the $x_1$-axis lies in $\operatorname{rad}(q_0)\cap\operatorname{rad}(q_1)$. Then $q_0+tq_1$ is degenerate, contrary to the hypothesis. Finally, if $\ell_{00}=\ell_{01}=0$, the differential of $q_0$ vanishes along $L$, again contradicting the hypothesis on $L$. Thus $\delta$ is nonzero. Let $U=\{\delta\ne0\}\subset\PP_K(E/A_2)$. For $[y]\in U$, the coordinates $x_0,x_1$ of a point of $X$ over $[y]$ satisfy
	$$
		M(y)\begin{pmatrix}x_0\\x_1\end{pmatrix}
		=-\begin{pmatrix}Q_0\\Q_1\end{pmatrix}.
	$$
	Since $\det M=\delta$, the matrix $M(y)$ is invertible for $[y]\in U$. Define
	$$
		\begin{aligned}
			\psi:U&\longrightarrow X\setminus L,\\
			[y]&\longmapsto
			[\ell_{01}Q_1-\ell_{11}Q_0:
			\ell_{10}Q_0-\ell_{00}Q_1:
			\delta y_1:\cdots:\delta y_{n-1}]
		\end{aligned}
	$$
	by solving the displayed linear system. The coordinates defining $\psi$ have degree three and do not vanish simultaneously on $U$, so $\psi$ is a morphism. Direct substitution gives
	$$
		\operatorname{pr}_L\circ\psi=\operatorname{id}_U,
		\qquad
		\psi\circ\operatorname{pr}_L
		=\operatorname{id}_{\operatorname{pr}_L^{-1}(U)}.
	$$
	Thus $\operatorname{pr}_L^{-1}(U)\simeq U$. Since $U$ is nonempty and $X$ is integral, $\operatorname{pr}_L^{-1}(U)$ is dense in $X$. Hence $\operatorname{pr}_L$ is birational.
\end{proof}

For a geometrically integral odd-degree curve, the following criterion is \cite[Lemma~3.20]{BF24}. The preceding argument gives the stated version for general cycle classes.

\begin{corollary}\label{cyc:field-rationality}
	Let $X\subset\PP_K(E)$ be an integral complete intersection of two quadrics, with $\dim_{K}E=n+1$ and $n\geq4$, and suppose that the generic member of its defining pencil is nondegenerate. Suppose also that the geometric nonsmooth locus of $X$ is finite. If $X$ has an $r$-cycle of odd degree with integer coefficients for some fixed $r\geq1$, then $X$ contains $\PP_K^r$ and is $K$-rational.
\end{corollary}
\begin{proof}
	Proposition~\ref{cyc:linear-space} gives $\PP_K^r\subset X$. Choose a $K$-line $L\subset\PP_K^r$. Since $L$ has dimension one and the geometric nonsmooth locus of $X$ is finite, $L$ is not contained in that locus. Proposition~\ref{cyc:projection} gives the birational map $X\dashrightarrow\PP_K^{n-2}$.
\end{proof}

The relative criterion follows by applying Corollary~\ref{cyc:field-rationality} over the function field of the base and using the projection formula to detect the degree of the restricted cycle.

Let $k$ be a field of characteristic different from $2$. For a dominant morphism $f:X\to B$ of $k$-varieties, put $K=k(B)$ and write $X_K=X\times_B\operatorname{Spec}K$ for its generic fiber.

\begin{theorem}\label{cyc:relative}
	Let $B$ be a projective $k$-variety of dimension $d$, and let $f:X\longrightarrow B$ be a dominant projective morphism of $k$-varieties. Assume that $X_K\subset\PP_K^n$, with $n\geq4$, is an integral complete intersection of two quadrics, that its defining pencil has nondegenerate generic member, and that its geometric nonsmooth locus is finite. Let $\mathcal H$ be a line bundle on $X$ such that $\mathcal H|_{X_K}\simeq\mathcal O_{X_K}(1)$. Fix $r\geq1$. Suppose that there exist a class $z\in\operatorname{CH}_{d+r}(X)$ and line bundles $\mathcal A_1,\ldots,\mathcal A_d$ on $B$ for which
	\begin{equation}\label{cyc:odd-intersection}
		\deg_{k}\!\left(
		c_1(\mathcal H)^r
		\prod_{i=1}^d c_1(f^*\mathcal A_i)\cap z\right)
		\quad\text{is odd}.
	\end{equation}
	Then $X_K$ contains a $K$-linear space $\PP_K^r$, and there is a birational map over $B$
	$$
		X\dashrightarrow B\times_{k}\PP_{k}^{n-2}.
	$$
	In particular, if $B$ is $k$-rational, then $X$ is $k$-rational.
\end{theorem}

\begin{proof}
	The restriction map $\operatorname{CH}_{d+r}(X)\to\operatorname{CH}_r(X_K)$ to the generic fiber sends $z$ to a class $z_K$. Define the integer
	$$
		\delta=\deg_K\bigl(c_1(\mathcal O_{X_K}(1))^r\cap z_K\bigr).
	$$
	We first prove the equality
	\begin{equation}\label{cyc:pushforward}
		f_*\bigl(c_1(\mathcal H)^r\cap z\bigr)=\delta[B]\quad\text{in }\operatorname{CH}_d(B).
	\end{equation}
	Choose a cycle $\sum_i m_i[W_i]$ representing $c_1(\mathcal H)^r\cap z$, where each $W_i\subset X$ is integral of dimension $d$. If $W_i$ dominates $B$, then $f|_{W_i}:W_i\to B$ is generically finite, since $\dim W_i=\dim B=d$. By the definition of proper pushforward,
	$$
		f_*[W_i]=[k(W_i):K]\,[B].
	$$
	The generic fiber of $f|_{W_i}$ is
	$$
		(W_i)_K=W_i\times_B\operatorname{Spec}K
		\simeq\operatorname{Spec}k(W_i),
	$$
	which is a closed point of $X_K$ of degree $[k(W_i):K]$ over $K$. If $W_i$ does not dominate $B$, then $f(W_i)$ is closed of dimension less than $d$ by properness, so $f_*[W_i]=0$. Since $f(W_i)$ does not contain the generic point of $B$, the fiber $(W_i)_K$ is empty. Since restriction to the generic fiber commutes with the action of first Chern classes, the cycle
	$$
		\sum_{f(W_i)=B}m_i[(W_i)_K]
	$$
	represents $c_1(\mathcal O_{X_K}(1))^r\cap z_K$. Therefore,
	$$
		\delta =\deg_K\bigl(c_1(\mathcal O_{X_K}(1))^r\cap z_K\bigr) =\sum_{f(W_i)=B}m_i[k(W_i):K],
	$$
	and hence
	$$
		f_*\bigl(c_1(\mathcal H)^r\cap z\bigr) =\delta[B] \quad\text{in }\mathrm{CH}_d(B).
	$$
	Applying the projection formula to \eqref{cyc:pushforward} gives
	$$
		\begin{aligned}
			\deg_{k} &\left( c_1(\mathcal H)^r \prod_{i=1}^d c_1(f^* \mathcal A_i) \cap z\right)\\
			&= \deg_{k} \left(\prod_{i=1}^d c_1(\mathcal A_i) \cap f_*( c_1( \mathcal H)^r \cap z)\right)\\
			&=\delta\, \deg_{k}\!\left(\prod_{i=1}^d c_1(\mathcal A_i)\cap[B]\right).
		\end{aligned}
	$$
	The left-hand side is odd by \eqref{cyc:odd-intersection}, so $\delta$ is odd. Corollary~\ref{cyc:field-rationality} gives $\PP_K^r\subset X_K$ and a birational map $X_K\dashrightarrow\PP_K^{n-2}$. Hence there is an isomorphism
	$$
		k(X)\simeq K(t_1,\ldots,t_{n-2})
		=k(B\times_k\PP_k^{n-2})
	$$
	over $K=k(B)$. This isomorphism defines a birational map $X\dashrightarrow B\times_k\PP_k^{n-2}$ over $B$. If $B$ is $k$-rational, then $k(X)$ is purely transcendental over $k$.
\end{proof}

\begin{remark}\label{cyc:hypotheses}
	Neither $X$ nor $B$ is required to be smooth. All products in \eqref{cyc:odd-intersection} are actions of first Chern classes, so they make sense on these possibly singular varieties. The class $z$ may have negative coefficients and vertical components. The hypothesis also forces the intersection number $\deg_{k}(\prod_{i=1}^d c_1(\mathcal A_i)\cap[B])$ to be odd. When $d=0$, the empty product is understood to be one.
\end{remark}

\begin{corollary}\label{cyc:projective-base}
	In the setting of Theorem~\ref{cyc:relative}, suppose that $B=\PP_{k}^d$, $n=5$, and $r=1$, and put $h_B=c_1(f^*\mathcal O_{\PP^d}(1))$. If a class $z\in\operatorname{CH}_{d+1}(X)$ satisfies
	$$
		\deg_{k}\bigl(c_1(\mathcal H)h_B^d\cap z\bigr)
		\equiv1\pmod2,
	$$
	then $X$ is $k$-rational.
\end{corollary}
\begin{proof}
	Take $\mathcal A_i=\mathcal O_{\PP^d}(1)$ for every $i\in\{1,\ldots,d\}$ in Theorem~\ref{cyc:relative}.
\end{proof}

For $d=3$, Corollary~\ref{cyc:projective-base} gives the cycle-theoretic form of \cite[Theorem~3.21]{BF24} stated as Corollary~\ref{pes:cycle-criterion}. It allows arbitrary classes in $\operatorname{CH}_4(X)$ and omits the geometric integrality assumption on the generic curves. This follows from Proposition~\ref{cyc:linear-space} and the projection formula.

\section{Complete intersections of an arbitrary number of quadrics}\label{higher:section}
We combine the odd-degree criterion for a quadric with hyperbolic reduction to obtain a rationality criterion for complete intersections of $m$ quadrics containing a $K$-linear space $\PP_K^{m-2}$. The hyperbolic reduction construction is standard; see \cite[Lemmas~3.3 and 3.4]{KS18}. We state the resulting quadric fibration and identify its generic fiber below. The coordinate construction and explicit birational maps are given in Appendix~\ref{app:hyperbolic}.

\begin{proposition}\label{higher:reduction}
	Let $K$ be a field of characteristic different from $2$, let $m\geq2$, and let $n\geq3m-2$. Let
	$$
		q_i:K^{n+1}\longrightarrow K,
		\qquad 1\leq i\leq m,
	$$
	be quadratic forms. Suppose that
	$$
		Y=V(q_1,\ldots,q_m)\subset\PP_K^n
	$$
	is a smooth integral complete intersection of $m$ quadrics and contains a linear subspace $\PP(U)\simeq\PP_K^{m-2}$ defined over $K$. Suppose also that the generic member of the linear system generated by $q_1,\ldots,q_m$ is nondegenerate. Then the following two statements hold:

    \begin{enumerate}
        \item There is an integral projective $K$-variety $Z$ and a 
                flat projective morphism
                $$
                    \pi:Z\longrightarrow\PP_K^{m-1}
                $$
        	    such that $Y$ and $Z$ are $K$-birational.
        \item Let $\eta$ be the generic point of $\PP_K^{m-1}$. Choose a nonzero vector $(\lambda_1,\ldots,\lambda_m)\in K(\PP_K^{m-1})^m$ representing $\eta$, and put
            	$$
            		q_\eta=\sum_{i=1}^m\lambda_i(q_i)_{K(\PP_K^{m-1})}.
            	$$
            	The quadratic form
            	$$
            		q_{\mathrm{red}}:U_{K(\PP_K^{m-1})}^{\perp_{q_\eta}}/U_{K(\PP_K^{m-1})}\longrightarrow K(\PP_K^{m-1}),
            		\qquad [v]\longmapsto q_\eta(v),
            	$$
                is nondegenerate, and
            	$$
            		Z_\eta\simeq V(q_{\mathrm{red}})
            		\subset\PP_{K(\PP_K^{m-1})}\!\left(U_{K(\PP_K^{m-1})}^{\perp_{q_\eta}}/U_{K(\PP_K^{m-1})}\right),
            		\qquad
            		q_\eta\simeq\HH^{\perp(m-1)}\perp q_{\mathrm{red}}.
            	$$
            	Here $Z_\eta=Z\times_{\PP_K^{m-1}}\operatorname{Spec}K(\PP_K^{m-1})$, and the orthogonal complement is taken with respect to $b_{q_\eta}$. In particular, $Z_\eta$ is a smooth quadric of dimension $n-2m+1$.
    \end{enumerate}
\end{proposition}

\begin{proof}
	Lemma~\ref{higher:quadric-family} in Appendix~\ref{app:hyperbolic} constructs $Z$ and $\pi$. Its proof identifies the reduced form $\overline q_\eta$ with $q_{\mathrm{red}}$ and gives $Z_\eta\simeq V(q_{\mathrm{red}})$ and the orthogonal decomposition~\eqref{higher:hyperbolic-decomposition}. Lemma~\ref{higher:birational} gives the $K$-birational map between $Y$ and $Z$.
\end{proof}

\begin{theorem}\label{higher:rationality}
	Under the hypotheses of Proposition~\ref{higher:reduction}, suppose that there exist an integer $r\geq m-1$ and a class $z\in\operatorname{CH}_r(Y)$ such that
	$$
		\deg_K\bigl(c_1(\cO_Y(1))^r\cap z\bigr)
	$$
	is odd.  Then $Y$ is $K$-rational.
\end{theorem}

\begin{proof}
	Use $q_\eta$ and $q_{\mathrm{red}}$ from Proposition~\ref{higher:reduction}. Write $j:Y_{K(\PP_K^{m-1})}\hookrightarrow V(q_\eta)$ for the inclusion. The projection formula and invariance of degree under field extension give
	$$
		\deg_{K(\PP_K^{m-1})}\bigl(c_1(\cO_{V(q_\eta)}(1))^r\cap j_*z_{K(\PP_K^{m-1})}\bigr)
		=\deg_K\bigl(c_1(\cO_Y(1))^r\cap z\bigr).
	$$
	Thus $j_*z_{K(\PP_K^{m-1})}$ has odd degree. Writing $i_W$ for the Witt index, Lemma~\ref{cyc:quadric} gives
	$$
		i_W(q_\eta)\geq r+1\geq m.
	$$
	By Proposition~\ref{higher:reduction}, we have
	$$
		q_\eta\simeq\HH^{\perp(m-1)}\perp q_{\mathrm{red}},
		\qquad
		i_W(q_{\mathrm{red}})=i_W(q_\eta)-(m-1)\geq1.
	$$
	The generic fiber $Z_\eta\simeq V(q_{\mathrm{red}})$ of $\pi:Z\to\PP^{m-1}_K$ therefore has a $K(\PP_K^{m-1})$-point. Since $Z_\eta$ is a smooth quadric of positive dimension, projection from a $K(\PP_K^{m-1})$-point of $Z_\eta$ proves that $Z_\eta$ is $K(\PP_K^{m-1})$-rational. Proposition~\ref{higher:reduction} now gives
	$$
		K(Y)\simeq K(Z)\simeq
		K(t_1,\ldots,t_{m-1})(v_1,\ldots,v_{n-2m+1}),
	$$
	where the $t_i$ and $v_j$ are algebraically independent over $K$. Thus $Y$ is $K$-rational of dimension $n-m$.
\end{proof}

\begin{remark}\label{higher:comparison}
	\begin{enumerate}
		\item For $K=\mathbf{C}$, $m=3$, and $n=7$, the rationality criterion using a surface cycle of odd degree is proved in \cite[Theorem~1.7]{RS23Q}. The quadric surface fibration birational to $Y$ is the construction recalled in \cite[Proposition~6]{HPT18}.  Theorem~\ref{higher:rationality} applies the odd-degree argument over a nonclosed field and for arbitrary $m$ in the stated range, using the specified linear subspace to perform hyperbolic reduction.
		\item The dimension condition $r \geq m-1$ on the cycle cannot be dropped: smooth complex intersections of three quadrics in $\PP^7$ contain lines, whereas a very general such intersection is not stably rational \cite[Theorem~1(1) and Proposition~6]{HPT18}.
		\item Over $\CC$, the dimension bound $n\geq3m-2$ is forced by the required odd-degree cycle. Indeed, put $d=\dim Y=n-m$ and $H=c_1(\cO_Y(1))$. An odd-degree cycle has $0\leq r\leq d$. If $d<2r$ and $c=d-r$, then $2c<d$; integral weak Lefschetz gives $H^{2c}(Y,\ZZ)=\ZZ\,H^c$. Thus the cycle class of $z$ is $aH^c$ for some $a\in\ZZ$, and its degree is $a\int_Y H^d=a2^m$, which is even. Consequently an odd-degree $r$-cycle requires $2r\leq d$. Together with $r\geq m-1$, this implies $n\geq3m-2$.
	\end{enumerate}
\end{remark}

The following corollary gives a rationality criterion for families whose generic fibers are complete intersections of $m$ quadrics.
\begin{corollary}\label{higher:relative}
	Let $k$ be a field of characteristic different from $2$, let $B$ be a projective $k$-variety of dimension $d$, and let $f:X\to B$ be a dominant projective morphism of $k$-varieties. Set $K=k(B)$.  Suppose that the generic fiber $X_K := X \times_{B} \operatorname{Spec}K\subset\PP^n_K$ satisfies the hypotheses of Proposition~\ref{higher:reduction}, and that a line bundle $\mathcal H$ on $X$ restricts to $\cO_{X_K}(1)$. Fix an integer $r\geq m-1$. Suppose that there exist a class $z\in\operatorname{CH}_{d+r}(X)$ and line bundles $\mathcal A_1,\ldots,\mathcal A_d$ on $B$ for which the integer
	$$
		\deg_{k}\left(c_1(\mathcal H)^r
		\prod_{i=1}^d f^*c_1(\mathcal A_i)\cap z\right)
	$$
	is odd. Then there is a birational map over $B$
	$$
		X\dashrightarrow B\times_{k}\PP^{n-m}_{k}.
	$$
	In particular, $X$ is $k$-rational if $B$ is $k$-rational.
\end{corollary}

\begin{proof}
	Let $z_K\in\operatorname{CH}_r(X_K)$ be the restriction of $z\in\operatorname{CH}_{d+r}(X)$, and set
	$$
		\delta=\deg_K\bigl(c_1(\mathcal H_K)^r\cap z_K\bigr).
	$$
	Then proper pushforward along $f:X\to B$ gives $f_*(c_1(\mathcal H)^r\cap z)=\delta[B]$ in $\operatorname{CH}_d(B)$. Hence the projection formula yields
	$$
		\deg_k\left(c_1(\mathcal H)^r\prod_{i=1}^d f^*c_1(\mathcal A_i)\cap z\right)
		=\delta\,\deg_k\left(\prod_{i=1}^d c_1(\mathcal A_i)\cap[B]\right).
	$$
	The left-hand side is odd by hypothesis, so $\delta$ is odd. Theorem~\ref{higher:rationality} then shows that $X_K$ is $K$-rational. Hence there is an isomorphism over $K=k(B)$
	$$
		k(X)\simeq K(t_1,\ldots,t_{n-m})
		=k\bigl(B\times_k\PP_k^{n-m}\bigr).
	$$
	This isomorphism gives a birational map $X\dashrightarrow B\times_k\PP_k^{n-m}$ over $B$. If $B$ is $k$-rational, then $k(X)$ is purely transcendental over $k$.
\end{proof}

\section{Constructing odd-degree cycles of dimension one}
\label{cub:section}

For the application to the rationality of Peskine sixfolds in Theorem~\ref{intro:peskine}, we construct odd-degree cycles of dimension one as follows. We consider a pencil of quadrics in $\PP_K^5$ for which three distinct members over $K^{\mathrm{sep}}$ each contain a chosen linear subspace isomorphic to $\PP_{K^{\mathrm{sep}}}^3$. We assume that the set of these three pairs, each consisting of a quadric and its chosen linear subspace, is Galois-invariant. From these linear subspaces, we construct a cycle of dimension one and degree three on the base locus of the pencil, defined over $K$. This provides the odd-degree cycle required by the rationality criterion in Section~\ref{cyc:section}.

Let $E$ be a six-dimensional vector space over a field $K$ of characteristic different from two. Fix a two-dimensional subspace $P\subset\operatorname{Sym}^2 E^\vee$, write $\Pi=\PP_K(P)$ for its pencil, and let
$$
	X=V(q_0,q_1)\subset\PP_K(E),
$$
where $(q_0,q_1)$ is any basis of $P$. For a field extension $F/K$, we write $X_F = X \times_K \operatorname{Spec}F$ for the base change.

\begin{lemma}\label{cub:intersection}
	Let $F/K$ be a field extension, and let $Q,Q'$ be distinct $F$-members of $\Pi_F$. Suppose that four-dimensional subspaces $A_4,A'_4\subset E \otimes_K F$ satisfy $\PP_F(A_4)\subset Q$ and $\PP_F(A'_4)\subset Q'$. Then
	$$
		\dim_F(A_4\cap A'_4)\geq2,
		\qquad \PP_F(A_4\cap A'_4)\subset X_F.
	$$
\end{lemma}

\begin{proof}
	The dimension formula gives
	$$
		\dim_F(A_4\cap A'_4)=8-\dim_F(A_4+A'_4)\geq2.
	$$
	Choose equations $q,q'\in P\otimes_K F$ for $Q$ and $Q'$. Since $Q\ne Q'$, the forms $q,q'$ are linearly independent and span $P\otimes_K F$. Thus $X_F=V(q,q')=Q\cap Q'$, and
	$$
		\PP_F(A_4\cap A'_4)\subset Q\cap Q'=X_F.
	$$
\end{proof}

Recall that a finite \'etale scheme $T$ of degree $n$ over a field $K$ splits over a separable closure $K^{\mathrm{sep}}$:
$$
	T\times_K\operatorname{Spec}K^{\mathrm{sep}} \simeq \bigsqcup_{i=1}^n\operatorname{Spec}K^{\mathrm{sep}}.
$$
In particular, the set $T(K^{\mathrm{sep}})$ has $n$ elements.

\begin{theorem}\label{cub:cubic-descent}
	Suppose that $q_0+tq_1$ is nondegenerate over $K(t)$. Let $T$ be a finite \'etale $K$-scheme of degree $3$, and let
	$$
		\theta:T\longrightarrow\Pi,
		\qquad w:T\longrightarrow\operatorname{Gr}_K(4,E)
	$$
	be morphisms. Assume that $\theta$ is injective on geometric points and that, for every $\tau\in T(K^{\mathrm{sep}})$,
	$$
		\PP_{K^{\mathrm{sep}}}(A_{4,\tau})
		\subset Q_\tau,
	$$
	where $A_{4,\tau}$ is defined by $w(\tau)$ and $Q_\tau$ by $\theta(\tau)$. Then $X$ has a one-cycle with integer coefficients of degree three and contains a $K$-line. If $T=\operatorname{Spec} K'$ for a separable cubic extension $K'/K$, then $X_{K'}$ contains a $K'$-line whose pushforward under the finite projection $X_{K'}\to X$ has degree three on $X$. If, in addition, $X$ is an integral complete intersection with finite geometric nonsmooth locus, then $X$ is $K$-rational.
\end{theorem}

\begin{proof}
	Write $T=\coprod_a\operatorname{Spec}K'_a$ and $T(K^{\mathrm{sep}})=\{\tau_1,\tau_2,\tau_3\}$. The map
	$$
		T(K^{\mathrm{sep}})\longrightarrow
		\{\text{unordered pairs of distinct elements of }T(K^{\mathrm{sep}})\},
		\qquad
		\tau_i\longmapsto T(K^{\mathrm{sep}})\setminus\{\tau_i\}
	$$
	is a Galois-equivariant bijection. For each component $\operatorname{Spec}K'_a$, choose a point $\tau_i$ above it, and put $\{j,h\}=\{1,2,3\}\setminus\{i\}$. The point $\tau_i$ determines an embedding $K'_a\hookrightarrow K^{\mathrm{sep}}$. Since $\theta$ and $w$ are defined over $K$, the subgroup $\operatorname{Gal}(K^{\mathrm{sep}}/K'_a)$ preserves $A_{4,\tau_j}\cap A_{4,\tau_h}$. Galois descent gives a subspace $B_a\subset E\otimes_K K'_a$ such that
	$$
		B_a\otimes_{K'_a}K^{\mathrm{sep}}
		=A_{4,\tau_j}\cap A_{4,\tau_h}.
	$$
	Lemma~\ref{cub:intersection} gives $\dim_{K'_a}B_a\geq2$. Moreover, $q_0|_{B_a}=q_1|_{B_a}=0$, since these restrictions vanish after the faithful extension of scalars to $K^{\mathrm{sep}}$. We choose a two-dimensional subspace $A_{2,a}\subset B_a$, and set $L_a=\PP_{K'_a}(A_{2,a})\subset X_{K'_a}$. Let
	$$
		\nu_a:X_{K'_a}\longrightarrow X
	$$
	be the finite base-change projection, and put
	$$
		z=\sum_a(\nu_a)_*[L_a]\in\operatorname{CH}_1(X).
	$$
	Since $\nu_a^*\mathcal O_X(1)|_{L_a}=\mathcal O_{L_a}(1)$, the projection formula gives
	$$
		\begin{aligned}
			\deg_K\bigl(c_1(\mathcal O_X(1))\cap z\bigr)
			&=\sum_a [K'_a:K] \,
			\deg_{K'_a}\bigl(c_1(\mathcal O_{L_a}(1))\cap[L_a]\bigr)\\
			&=\sum_a[K'_a:K]=3.
		\end{aligned}
	$$
	Hence Proposition~\ref{cyc:linear-space}, with $r=1$, gives a $K$-line on $X$. If $T=\operatorname{Spec}K'$, there is a single component, so $L_a\subset X_{K'}$ and $z=(\nu_a)_*[L_a]$ have the required properties. If $X$ is an integral complete intersection with finite geometric nonsmooth locus, Corollary~\ref{cyc:field-rationality} shows that $X$ is $K$-rational.
\end{proof}

Let $q_\theta$ denote the quadratic form selected by the morphism $\theta:T\to\Pi$; it is well defined up to a unit. The condition that $q_\theta$ vanish identically on the universal rank four subbundle defines a closed subscheme
$$
	\mathcal F_T\subset
	\operatorname{Gr}_T(4,E\otimes_K\mathcal O_T).
$$
Thus $\mathcal F_T$ parametrizes projective spaces of dimension three on the corresponding members of the pencil.

\begin{corollary}\label{cub:conic}
	Let $\Pi=\mathbf P_K(P)$ be the pencil of quadrics generated by $q_0$ and $q_1$ considered in Theorem~\ref{cub:cubic-descent}, and let $\theta:T\to\Pi$ be as in that theorem. Suppose that there is a smooth proper family of conics
	$$
		\pi_C:C\longrightarrow T
	$$
	with geometrically connected fibers and a $T$-morphism
	$$
		\iota:C\longrightarrow\mathcal F_T.
	$$
	If $\pi_C$ has a section, then $X$ has a one-cycle of degree three and contains a $K$-line. If, in addition, $X$ is an integral complete intersection with finite geometric nonsmooth locus, then $X$ is $K$-rational. In particular, if $T=\operatorname{Spec}K'$ and $C$ has a $K'$-rational point, then $X$ has a one-cycle of degree three and contains a $K$-line.
\end{corollary}
\begin{proof}
	Let $s:T\to C$ be a section of $\pi_C$, and define $w$ by the composition
	$$
		w:T\xrightarrow{s}C\xrightarrow{\iota}\mathcal F_T
		\longrightarrow\operatorname{Gr}_K(4,E).
	$$
	Since $\iota\circ s$ is a section of $\mathcal F_T\to T$, the morphisms $\theta$ and $w$ satisfy the hypotheses of Theorem~\ref{cub:cubic-descent}. The assertions follow from that theorem.
\end{proof}

The section hypothesis in Corollary~\ref{cub:conic} requires a rational point on the conic over the residue field of each closed point of $T$. The proof uses only the resulting section of $\mathcal F_T\to T$. In Section~\ref{pes:section}, the conic is the pullback of the Brauer--Severi conic $\pi:D\to S$ along a morphism $\operatorname{Spec}\CC(\Gamma)\to S$ from the generic point of the incidence variety $\Gamma$, and the vanishing of its Brauer class provides this section. The proof of Theorem~\ref{pes:beta-rational} applies the argument of Theorem~\ref{cub:cubic-descent} directly, without using Corollary~\ref{cub:conic} and Corollary~\ref{cub:relative-cycle}.

\begin{corollary}\label{cub:relative-cycle}
	Let $f:X\to\PP_{k}^d$ and $\mathcal H$ satisfy the hypotheses of Theorem~\ref{cyc:relative} with $n=5$, and put $K=k(\PP_{k}^d)$. If the pencil defining $X_K$ satisfies the hypotheses of Theorem~\ref{cub:cubic-descent}, or those of Corollary~\ref{cub:conic} with $\pi_C:C\to T$ admitting a section, then there is a subvariety $Z_L\subset X$ of dimension $d+1$ satisfying
	$$
		\deg_{k}\bigl(
		c_1(\mathcal H)c_1(f^*\mathcal O_{\PP^d}(1))^d\cap[Z_L]
		\bigr)=1.
	$$
	Consequently, $X$ is $k$-rational.
\end{corollary}

\begin{proof}
	Theorem~\ref{cub:cubic-descent} or Corollary~\ref{cub:conic} gives a line $L\subset X_K$. Let $Z_L$ be the reduced closure of $L$ in $X$. Then $Z_L$ is integral of dimension $d+1$ and has generic fiber $L$. Since $\deg_K\mathcal O_L(1)=1$, formula~\eqref{cyc:pushforward} with $r=1$ gives
	$$
		f_*\bigl(c_1(\mathcal H)\cap[Z_L]\bigr)=[\PP_{k}^d].
	$$
	The projection formula now gives
	$$
		\begin{aligned}
			\deg_k\bigl(c_1(\mathcal H)c_1(f^*\mathcal O_{\PP^d}(1))^d\cap[Z_L]\bigr)
			&=\deg_k\bigl(c_1(\mathcal O_{\PP^d}(1))^d\cap[\PP_k^d]\bigr)\\
			&=1.
		\end{aligned}
	$$
	Corollary~\ref{cyc:projective-base} therefore shows that $X$ is $k$-rational.
\end{proof}

\section{Application to Peskine sixfolds}\label{pes:section}
Under Assumption~\ref{pes:generic}, we show that the Peskine sixfold $X_1^\sigma$ is rational whenever the Brauer class $\beta$ on the associated K3 surface vanishes. The K3 surface, the conic, and the birational model of the Peskine sixfold are the constructions of Benedetti--Song and Benedetti--Faenzi recalled below. The additional construction uses an incidence variety of degree $3$ to turn a rational section of the conic into a line on the generic fiber after a cubic field extension.

\subsection{The setup}\label{pes:geometric-setup}

Let $\sigma\in\mathcal D^{1,6,10}$ be general (see Assumption~\ref{pes:generic} below), and let $X_1^\sigma$ and $X_6^\sigma$ be as in \eqref{def:Peskine} and \eqref{def:Debarre-Voisin}. We give $X_1^\sigma$ its Pfaffian scheme structure and $X_6^\sigma$ its zero scheme structure. We fix a flag $V_1\subset V_6\subset V_{10}$ satisfying
\begin{equation}\label{pes:flag}
	\dim V_i=i\quad(i=1,6),\qquad \sigma(V_1,V_6,V_{10})=0.
\end{equation}
The divisor $\mathcal D^{1,6,10}$ parametrizes trivectors admitting such a flag. For general $\sigma$ the flag is unique, as recalled in \cite[\S3.1]{BF24}. We choose a vector $0\ne v_1\in V_1$. For the quotient map $\rho_6:V_{10}\to V_{10}/V_6$, equation~\eqref{pes:flag} gives the alternating form
$$
	\omega:(V_{10}/V_6)\times(V_{10}/V_6)\longrightarrow\CC,
	\qquad \omega(\rho_6(x),\rho_6(y))=\sigma(v_1,x,y).
$$
Its definition is independent of the lifts, since contraction with $v_1$ vanishes whenever one argument belongs to $V_6$. Rescaling $v_1$ rescales $\omega$ and leaves its isotropic subspaces unchanged. The uniqueness of the distinguished flag implies that the kernel of $\sigma(v_1,-,-)$ is exactly $V_6$, so $\omega$ is nondegenerate. Consequently
$$
	Q=\operatorname{LGr}(2,V_{10}/V_6,\omega)
	\subset\operatorname{Gr}(2,V_{10}/V_6)
$$
is a smooth quadric threefold in its Pl\"ucker embedding. A point $[U_8]\in Q$ denotes the Lagrangian plane $U_8/V_6$, where $V_6\subset U_8\subset V_{10}$ and $\sigma(V_1,U_8,U_8)=0$.

On $Q$, let $\mathcal U_8\subset V_{10}\otimes\cO_Q$ be the inverse image of the tautological bundle $\mathcal U_{8/6}\subset(V_{10}/V_6)\otimes\cO_Q$ under the quotient map
$$
	\rho_6\otimes\cO_Q:V_{10}\otimes\cO_Q
	\longrightarrow(V_{10}/V_6)\otimes\cO_Q.
$$
The restrictions of $\sigma$ descend to a section of $\bigwedge^3(\mathcal U_8/V_1)^\vee$. The degree seven relative invariant of three-forms in dimension seven, recalled in \cite[Appendix~C]{BS21}, takes values in the third power of the dual determinant. It therefore gives, up to a fixed nonzero scalar, a section
\begin{equation} \label{k3:section}
	c_\sigma\in H^0\!\left(Q, \det(V_6/V_1)^{\vee\otimes3}\otimes\cO_Q(3)\right).
\end{equation}
Here we used $\det(\mathcal U_8/V_1)^\vee \simeq\det(V_6/V_1)^\vee\otimes\det\mathcal U_{8/6}^\vee$. After fixing the constant determinant factor, write $c_\sigma\in H^0(Q,\cO_Q(3))$ and $S=V_Q(c_\sigma)$. Now we define the set
$$
	D=\{[U_6]\in X_6^\sigma\mid V_1\subset U_6\}.
$$
For general $\sigma$, the section $c_\sigma$ is nonzero, $S$ is a smooth K3 surface of degree six, and
$$
	\pi:D\longrightarrow S,\qquad
	[U_6]\longmapsto[(U_6+V_6)/V_6]
$$
is a Brauer--Severi conic, and $\dim(U_6\cap V_6)=4$ for every geometric point $[U_6]\in D$; see \cite[Propositions~5.5, 5.7 and 5.8]{BS21} and \cite[Propositions~3.1 and 3.2]{BF24}. In \cite{BF24}, the surface $S$ is defined by an incidence condition in $\operatorname{Gr}(2,V_{10}/V_6)$. Its identification with the zero scheme $V_Q(c_\sigma)$ is given in the proof of \cite[Proposition~5.5]{BS21}. Denote the Brauer class of $\pi$ by $\beta\in\operatorname{Br}(S)[2]$. If $\beta=0$, its generic conic has a rational point, and hence $\pi$ has a rational section $r:S\dashrightarrow D$ with $\pi\circ r=\operatorname{id}_S$. The dashed arrow $r$ in the following diagram exists only when $\beta=0$.
\begin{equation}
	\begin{tikzcd}[column sep=large, row sep=large]
		D \arrow[r, hook] \arrow[d, shift right=1.2ex, "\pi"'] &
		X_6^\sigma \arrow[r, hook] &
		\operatorname{Gr}(6,V_{10})
		\\
		S \arrow[r, hook] \arrow[u, dashed, shift right=1.2ex, "r\;(\beta=0)"'] &
		Q \arrow[r, hook] &
		\operatorname{Gr}(2,V_{10}/V_6)
	\end{tikzcd}
\end{equation}

The base of the Peskine fibration is
$$
	B=\PP(V_{10}/V_6),\qquad K=\CC(B).
$$
On $B$ let $\cO_B(-1)\subset(V_{10}/V_6)\otimes\cO_B$ be the tautological line bundle. Write $\cO_B(-1)^\perp$ for the kernel of the map
$$
	(V_{10}/V_6)\otimes\cO_B\longrightarrow\cO_B(1),
	\qquad x\longmapsto\omega(x,-)|_{\cO_B(-1)}.
$$
Following \cite[Remark~3.7 and \S3.2]{BF24}, let $\mathcal V_7$ and $\mathcal V_9$ be the inverse images of $\cO_B(-1)$ and $\cO_B(-1)^\perp$, respectively, under $\rho_6\otimes\cO_B$. These bundles have ranks seven and nine. The rank-two null-correlation bundle is
$$
	\mathcal N=\mathcal V_9/\mathcal V_7
	\simeq\cO_B(-1)^\perp/\cO_B(-1).
$$
In the quotient $\mathcal V_7/V_1$, we write $V_1$ for the constant subbundle $V_1\otimes\cO_B$. Each point $b\in B$ is represented by a line in $V_{10}/V_6$. Let $U_7\subset V_{10}$ be its inverse image, so $b=[U_7/V_6]$. As in \cite[\S3.2]{BF24}, we also write $b=[U_7]$. Define
$$
	U_7^\perp=\rho_6^{-1}\bigl((U_7/V_6)^{\perp_\omega}\bigr).
$$
The fibers of $\mathcal V_7$, $\mathcal V_9$, $\mathcal V_7/V_1$, and $\mathcal N$ at $b$ are $U_7$, $U_7^\perp$, $U_7/V_1$, and $\mathcal N_b=U_7^\perp/U_7$, respectively. In particular, $\PP(\mathcal N_b)$ parametrizes the eight-dimensional subspaces $U_8$ with $U_7\subset U_8\subset U_7^\perp$. The point corresponding to $U_8$ is $[U_8/U_7]$. For a morphism $\Phi$ of vector bundles on $B$, we write $\Phi_b$ for the induced map on the fibers at $b$.

We define the birational model using the projection
$$
	p_B:\PP_B(\mathcal V_7/V_1)\longrightarrow B.
$$
On $\PP_B(\mathcal V_7/V_1)$ let $\mathcal U_1=\cO_{\PP_B(\mathcal V_7/V_1)}(-1)$ be the tautological line subbundle of $p_B^*(\mathcal V_7/V_1)$. Let $\mathcal V_2\subset p_B^*\mathcal V_7$ be its inverse image under the quotient map. At the point $[U_2/V_1]\in\PP(U_7/V_1)$, the fibers of $\mathcal U_1$ and $\mathcal V_2$ are $U_2/V_1$ and $U_2$, respectively.
The $\PP^1$-bundle $\PP_{\PP_B(\mathcal V_7/V_1)}(\mathcal V_2)\to\PP_B(\mathcal V_7/V_1)$ comes with the morphisms
$$
	\begin{aligned}
		p_{1,2}:\PP_{\PP_B(\mathcal V_7/V_1)}(\mathcal V_2)
		&\longrightarrow\PP_B(\mathcal V_7/V_1),\\
		\rho:\PP_{\PP_B(\mathcal V_7/V_1)}(\mathcal V_2)
		&\longrightarrow\PP(V_{10}),
		\qquad [\ell\subset U_2\subset U_7]\longmapsto[\ell].
	\end{aligned}
$$

\begin{definition}[{\cite[Definitions~3.9 and 3.12]{BF24}}]\label{def:birational-models}
	Define $X_{1.1}^\sigma$ to be the closure in $\PP_{\PP_B(\mathcal V_7/V_1)}(\mathcal V_2)$ of the inverse image of $X_1^\sigma\setminus\mathbf P(V_6)$ under $\rho$. Let $X_{1.2}^\sigma$ be the reduced image of $X_{1.1}^\sigma$ under $p_{1,2}$.
\end{definition}

Recall that $K=\mathbf C(B)$ is the function field of $B$. We write $\varpi$ for the restriction of $p_B$ to $X_{1.2}^\sigma$ and $X_{1.2/K}^\sigma$ for the generic fiber of $\varpi$:
$$
    \varpi:X_{1.2}^\sigma\longrightarrow B, \qquad
    X_{1.2/K}^\sigma=X_{1.2}^\sigma\times_B\operatorname{Spec}K \subset\mathbf P_K\bigl((\mathcal V_7/V_1)_K\bigr).
$$
For a field extension $F/K$, write $X_{1.2/F}^\sigma=X_{1.2/K}^\sigma\times_K F$. Here $(\mathcal V_7)_K$ is the fiber of $\mathcal V_7$ at the generic point of $B$, and $(\mathcal V_7/V_1)_K=(\mathcal V_7)_K/(V_1\otimes_{\CC}K)$ has dimension six over $K$. By \cite[Propositions~3.10 and 3.17--3.19]{BF24}, the variety $X_{1.2}^\sigma$ is integral and birational to $X_1^\sigma$. Its generic fiber $X_{1.2/K}^\sigma$ is a three-dimensional complete intersection of the two quadrics in the Pfaffian pencil described in Lemma~\ref{pes:pfaffian}. The geometric singular locus $\operatorname{Sing}(X_{1.2/\bar K}^\sigma)$ is finite, and the generic member of the pencil is smooth; equivalently, its determinant polynomial is nonzero.

The genericity conditions used below are the following.

\begin{assumption}\label{pes:generic}
	Fix a nonempty Zariski open subset $U\subset\mathcal D^{1,6,10}$ on which the flag~\eqref{pes:flag} is unique and the conclusions of \cite[Propositions~3.1, 3.2, 3.10, and 3.17--3.19]{BF24} hold. Throughout this section, we assume that $\sigma\in U$.
\end{assumption}

Such a nonempty Zariski open subset $U$ exists by \cite[\S3.1 and Propositions~3.1, 3.2, 3.10, and 3.17--3.19]{BF24}. We say that $\sigma$ is \emph{general as in Section~\ref{pes:geometric-setup}} if $\sigma\in U$.

\begin{remark}\label{rem:generic}
	Under Assumption~\ref{pes:generic}, the following properties hold:
	\begin{enumerate}
		\item the flag~\eqref{pes:flag} is unique;
		\item $c_\sigma\ne0$, the surface $S$ is a smooth K3 surface of degree six, and $\pi:D\to S$ is a Brauer--Severi conic;
		\item $X_{1.2}^\sigma$ is integral and birational to $X_1^\sigma$;
		\item $X_{1.2/K}^\sigma$ is a complete intersection of two quadrics in $\PP_K\bigl((\mathcal V_7/V_1)_K\bigr)$ whose geometric nonsmooth locus is finite and whose defining pencil has a nondegenerate generic member.
	\end{enumerate}
\end{remark}

We summarize the varieties and maps introduced above in the following diagrams.
\begin{equation}
	\begin{tikzcd}[column sep=large, row sep=large]
		X_1^\sigma \arrow[d, hook] &
		X_{1.1}^\sigma \arrow[l, "\rho"'] \arrow[r, "p_{1,2}"] \arrow[d, hook] &
		X_{1.2}^\sigma \arrow[rd, "\varpi"] \arrow[d, hook]
		\\
		\mathbf P(V_{10}) &
		\mathbf P_{\mathbf P_B(\mathcal V_7/V_1)}(\mathcal V_2) \arrow[l, "\rho"'] \arrow[r, "p_{1,2}"'] &
		\mathbf P_B(\mathcal V_7/V_1) \arrow[r, "p_B"'] &
		B
	\end{tikzcd}
\end{equation}
The arrows labelled $\rho$ and $p_{1,2}$ denote their restrictions to $X_{1.1}^\sigma$.
\begin{equation}
	\begin{tikzcd}[column sep=large, row sep=large]
		X_{1.2/K}^\sigma \arrow[r, hook] \arrow[d] & \mathbf P_K\bigl((\mathcal V_7/V_1)_K\bigr) \arrow[r] \arrow[d] &
		\operatorname{Spec}K \arrow[d] \\
		X_{1.2}^\sigma \arrow[r, hook] \arrow[rr, bend right=25, "\varpi"'] &
		\mathbf P_B(\mathcal V_7/V_1) \arrow[r, "p_B"] & B
	\end{tikzcd}
\end{equation}
Both squares in the last diagram are Cartesian.

We first express the Pfaffian pencil described in \cite[Proposition~3.18]{BF24} and in the discussion preceding \cite[Proposition~3.19]{BF24} as a morphism of vector bundles. We then construct an incidence variety $\Gamma$ and show that its projection $\lambda:\Gamma\to B$ is generically finite of degree three. We also construct a morphism from the pulled-back conic to the Fano scheme of projective three-spaces.

\begin{lemma}\label{pes:pfaffian}
	The trivector $\sigma \in U$ determines a morphism of vector bundles
	$$
		\Phi:\mathcal N\longrightarrow
		\operatorname{Sym}^2(\mathcal V_7/V_1)^\vee\otimes\det(\mathcal V_7/V_1)^\vee.
	$$
	For fixed $b=[U_7]\in B$ and $U_7\subset U_8\subset U_7^\perp$, the restriction of $\Phi_b$ to $\ell_{U_8}=U_8/U_7$ is the quadratic Pfaffian section
	\begin{equation}  \label{pes:quadratic-pfaffian}
		s_{U_8}\in H^0\!\left(\PP(U_7/V_1),\cO(2)\otimes
		\det (U_7/V_1)^\vee\otimes\ell_{U_8}^\vee\right).
	\end{equation}
	At the generic point of $B$, write $\mathcal N_K$ for the fiber of $\mathcal N$. The induced map $\Phi_K:\mathcal N_K\to\operatorname{Sym}^2(\mathcal V_7/V_1)_K^\vee\otimes\det (\mathcal V_7/V_1)_K^\vee$ has rank two and identifies $\Pi_K:=\PP_K(\mathcal N_K)$ with the pencil defining $X_{1.2/K}^\sigma$.
\end{lemma}

\begin{proof}
	Fix $b=[U_7]\in B$ and $U_7\subset U_8\subset U_7^\perp$. The two-dimensional subspace $U_8/V_6$ contains $U_7/V_6$ and is contained in $(U_7/V_6)^{\perp_\omega}$, so it is Lagrangian. Thus $\sigma(V_1,U_8,U_8)=0$. At $[U_2/V_1]\in\PP(U_7/V_1)$ there is a canonical alternating map with values in the line $(U_2/V_1)^\vee$,
	$$
		\widetilde\alpha_{U_2,U_8}:\bigwedge^2(U_8/U_2)
		\longrightarrow(U_2/V_1)^\vee,
	$$
	given by
	$$
		(x+U_2)\wedge(y+U_2)\longmapsto
		\bigl(u+V_1\longmapsto\sigma(u,x,y)\bigr).
	$$
	Since $\sigma(V_1,U_8,U_8)=0$ and $\dim(U_2/V_1)=1$, we also have $\sigma(U_2,U_2,U_8)=0$. These two equalities show that $\widetilde\alpha_{U_2,U_8}$ is independent of the lifts $u,x,y$. Choosing a nonzero vector $\bar u=u+V_1\in U_2/V_1$ defines an isomorphism
	$$
		\operatorname{ev}_{\bar u}:(U_2/V_1)^\vee\longrightarrow\mathbf C,
		\qquad
		\lambda\longmapsto\lambda(\bar u).
	$$
	The composite $\operatorname{ev}_{\bar u}\circ\widetilde\alpha_{U_2,U_8}$ is the scalar-valued alternating form on $U_8/U_2$ given by $(\bar x,\bar y)\mapsto\sigma(u,x,y)$. For $c\in\CC^\times$, we have $\operatorname{ev}_{c\bar u}\circ\widetilde\alpha_{U_2,U_8}=c(\operatorname{ev}_{\bar u}\circ\widetilde\alpha_{U_2,U_8})$; the line-valued form $\widetilde\alpha_{U_2,U_8}$ does not depend on $\bar u$. Put $W=U_8/V_1$, as in the proof of \cite[Proposition~3.18]{BF24}. On $\PP(U_7/V_1)$ use $\mathcal U_1$ for the restriction of the relative tautological line bundle. Define $\mathcal W_6$ by the exact sequence
	$$
		0 \longrightarrow \mathcal U_1 \longrightarrow W\otimes\cO \longrightarrow \mathcal W_6 \longrightarrow 0.
	$$
	The fiber of $\mathcal W_6$ at $[U_2/V_1]$ is $W/(U_2/V_1)=U_8/U_2$. As $[U_2/V_1]$ varies in $\mathbf P(U_7/V_1)$, the line-valued alternating forms $\widetilde\alpha_{U_2,U_8}$ define a morphism of vector bundles
	$$
		\widetilde\alpha:\bigwedge^2\mathcal W_6\longrightarrow \mathcal U_1^\vee,
	$$
	whose fiber at $[U_2/V_1]$ is $\widetilde\alpha_{U_2,U_8}$. The Pfaffian of $\widetilde\alpha$ is the section
	$$
		\mathrm{Pf}(\widetilde\alpha) = \frac{1}{3!} \widetilde\alpha^{\wedge 3}
	$$
	of $\det\mathcal W_6^\vee\otimes(\mathcal U_1^\vee)^{\otimes3}$. The sequence defining $\mathcal W_6$ and the exact sequence
	$$
		0\longrightarrow U_7/V_1\longrightarrow U_8/V_1
		\longrightarrow\ell_{U_8}\longrightarrow0
	$$
	give
	$$
		\begin{aligned}
			\det\mathcal W_6^\vee\otimes(\mathcal U_1^\vee)^{\otimes3}
			&\simeq\det(U_7/V_1)^\vee\otimes\ell_{U_8}^\vee
			\otimes \mathcal U_1\otimes(\mathcal U_1^\vee)^{\otimes3}\\
			&\simeq\cO(2)\otimes\det(U_7/V_1)^\vee\otimes\ell_{U_8}^\vee.
		\end{aligned}
	$$
	In the notation of \cite[Lemma~3.15]{BF24}, $\widetilde\alpha_{U_2,U_8}$ is the restriction of $\sigma''([U_2\subset U_7])$ to $\bigwedge^2(U_8/U_2)$. Indeed, the ambiguity in a representative of $\sigma''([U_2\subset U_7])$ lies in $(U_2/V_1)^\vee\otimes\bigwedge^2(U_7^\perp/U_7)^\vee$, whose restriction to the line $U_8/U_7$ is zero. Define $s_{U_8}=\mathrm{Pf}(\widetilde\alpha)$ under the determinant identification above. Then $s_{U_8}$ is the section in~\eqref{pes:quadratic-pfaffian}. At a point $[U_2/V_1]$, the Pfaffian $\mathrm{Pf}(\widetilde\alpha)$ vanishes if and only if the alternating form $\widetilde\alpha_{U_2,U_8}$ has rank at most four. To describe the dependence of $s_{U_8}$ on $U_8/U_7$, let $\mathcal W_7$ be the inverse image of $\cO_{\PP(\mathcal N_b)}(-1)$ under the quotient map
	$$
		(U_7^\perp/V_1)\otimes\cO_{\PP(\mathcal N_b)}
		\longrightarrow \mathcal N_b\otimes\cO_{\PP(\mathcal N_b)}.
	$$
	The bundle $\mathcal W_7$ has fiber $U_8/V_1$ at $[U_8/U_7]$ and fits into the exact sequence
	$$
		0\longrightarrow(U_7/V_1)\otimes\cO
		\longrightarrow\mathcal W_7
		\longrightarrow\cO_{\PP(\mathcal N_b)}(-1)\longrightarrow0.
	$$
	Applying the same Pfaffian construction on $\PP(U_7/V_1)\times\PP(\mathcal N_b)$, with $U_8/V_1$ replaced by $\mathcal W_7$, gives a section
	$$
		s_b\in H^0\!\left(\PP(U_7/V_1)\times\PP(\mathcal N_b),
		\cO(2)\boxtimes\cO(1)\otimes\det(U_7/V_1)^\vee\right).
	$$
	The section $s_b$ restricts to $s_{U_8}$ over $[U_8/U_7]$. The natural isomorphism
	$$
		\begin{aligned}
			&H^0\!\left(\PP(U_7/V_1)\times\PP(\mathcal N_b),
			\cO(2)\boxtimes\cO(1)\otimes\det(U_7/V_1)^\vee\right)\\
			&\qquad\simeq\operatorname{Hom}\!\left(\mathcal N_b,
			\operatorname{Sym}^2(U_7/V_1)^\vee\otimes\det(U_7/V_1)^\vee\right)
		\end{aligned}
	$$
	associates a linear map $\Phi_b$ to $s_b$. Contraction by $\sigma$, passage to the quotient bundle, and taking the exterior cube commute with changes of local frame. Hence the maps $\Phi_b$ are the fibers of the bundle morphism $\Phi$. We next show that $X_{1.2/K}^\sigma=V(\operatorname{Im}\Phi_K)$ as schemes. For each line $\ell_{U_8}=U_8/(\mathcal V_7)_K\subset \mathcal N_K$, the restriction $\Phi_K|_{\ell_{U_8}}$ gives the quadratic Pfaffian section $s_{U_8}$. By the proof of \cite[Proposition~3.18]{BF24} and the argument preceding \cite[Proposition~3.19]{BF24}, the common zero scheme of the Pfaffian sections $s_{U_8}$ is $X_{1.2/K}^\sigma$. Thus
	$$
		X_{1.2/K}^\sigma=V(\operatorname{Im}\Phi_K)
		\subset\mathbf P_K\bigl((\mathcal V_7/V_1)_K\bigr)
	$$
	as schemes. Since $X_{1.2/K}^\sigma$ has codimension two and $\dim_K \mathcal N_K=2$, the map $\Phi_K$ has rank two.
\end{proof}

For a nonzero $t\in \mathcal N_K$, choosing a trivialization of $\det (\mathcal V_7/V_1)_K^\vee$ identifies $\Phi_K(t)$ with a quadratic form $q_t$ on $(\mathcal V_7/V_1)_K$; its zero scheme $Q_t$ depends only on $[t]\in\Pi_K$. Two distinct parameters give independent equations, so
\begin{equation}\label{pes:base-scheme}
	Q_t\cap Q_{t'}=X_{1.2/K}^\sigma\qquad([t]\ne[t'])
\end{equation}
as schemes. Using the restriction $\mathcal U_{8/6}|_S$ of the tautological rank-two bundle, form the incidence variety
$$
	\Gamma=\PP_S(\mathcal U_{8/6}|_S)
	=\{(s=[U_8],b=[U_7])\mid U_7\subset U_8\}.
$$
Its two projections are $\pi_S:\Gamma\to S$, $(s,b)\mapsto s$, and $\lambda:\Gamma\to B$, $(s,b)\mapsto b$.

\begin{lemma}\label{pes:incidence-variety}
	The variety $\Gamma$ is integral and $\lambda:\Gamma\to B$ is generically finite of degree three. Its generic fiber is $\operatorname{Spec}\CC(\Gamma)$, and $\CC(\Gamma)/K$ is a separable cubic extension. Let $b$ be the generic point of $B$, and put $\Gamma_K=\Gamma\times_B\operatorname{Spec}K$. The morphism
    $$
    	\Gamma_K\longrightarrow\Pi_K,
    	\qquad (s=[U_8],b)\longmapsto[U_8/(\mathcal V_7)_K]
    $$
    is a closed immersion. The image of $\Gamma_K$ is reduced and has degree three over $K$.
\end{lemma}

\begin{proof}
	Let
	$$
		\mathcal I=\{([U_8],[U_7])\in Q\times B\mid U_7\subset U_8\},
		\qquad p:\mathcal I\to B,\quad q:\mathcal I\to Q.
	$$
	The fibers of the projections are
	$$
		p^{-1}([U_7])=\PP(U_7^\perp/U_7),
		\qquad q^{-1}([U_8])=\PP(U_8/V_6).
	$$
	Thus $p$ and $q$ are $\PP^1$-bundles, and $\mathcal I$ is integral. The restriction of $q$ to each fiber of $p$ embeds that fiber as a Pl\"ucker line in $Q$. Hence, for $\mathcal I_K=\mathcal I\times_B\operatorname{Spec}K\simeq\PP_K^1$, the section $q^*c_\sigma$ restricts to
	$$
		(q^*c_\sigma)_K\in H^0(\mathcal I_K,\cO_{\mathcal I_K}(3)).
	$$
	Since $c_\sigma\ne0$ and $q$ is surjective, $q^*c_\sigma\ne0$. The restriction $(q^*c_\sigma)_K$ is also nonzero, because $\mathcal I$ is integral and its generic point lies in $\mathcal I_K$. Therefore
	$$
		\Gamma_K=V_{\mathcal I_K}((q^*c_\sigma)_K)
	$$
	has length three over $K$. The variety $\Gamma$ is integral, since it is a projective bundle over the integral surface $S$. Thus $\Gamma_K=\operatorname{Spec}\CC(\Gamma)$ and $[\CC(\Gamma):K]=3$. The extension $\CC(\Gamma)/K$ is separable because $\operatorname{char}K=0$. Finally, the quotient $\mathcal N_b=U_7^\perp/U_7$ is naturally isomorphic to $(U_7/V_6)^{\perp_\omega}/(U_7/V_6)$. These isomorphisms identify $\mathcal I_K$ with $\Pi_K$ and embed $\Gamma_K$ as a reduced closed subscheme of degree three.
\end{proof}

\begin{lemma}\label{pes:conic-spaces}
	Let $(s=[U_8],b=[U_7])$ be a point of $\Gamma$ with values in a field extension of $\CC$, and let $[U_6]\in D_s$. Then
	$$
		A_4=(U_6\cap U_7)/V_1\subset U_7/V_1
	$$
	has dimension four and $\PP(A_4)\subset Q_{U_8}$, where $Q_{U_8}$ is the quadric defined by~\eqref{pes:quadratic-pfaffian}. In particular, over the generic incidence point there is a morphism
	\begin{equation}\label{pes:conic-to-fano}
		D\times_S\Spec \CC(\Gamma)\longrightarrow
		\mathcal{M}_3(Q_{U_8}/\CC(\Gamma)),
		\qquad [U_6]\longmapsto[\PP((U_6\cap U_7)/V_1)].
	\end{equation}
	Here we use the Fano scheme of projective three-spaces defined in Section~\ref{pre:fano}.
\end{lemma}

\begin{proof}
	By the definition of $\pi$, we have $U_6+V_6=U_8$. If $U_6\subset U_7$, then $U_8=U_6+V_6\subset U_7$, contrary to $\dim U_8=8>7=\dim U_7$. Since $U_7$ is a hyperplane in $U_8$, it follows that
	$$
		\dim(U_6\cap U_7)=5,
		\qquad\dim A_4=\dim((U_6\cap U_7)/V_1)=4.
	$$
	For each $V_1\subset U_2\subset U_6\cap U_7$, the alternating form $\widetilde\alpha_{U_2,U_8}$ vanishes on the four-dimensional subspace $U_6/U_2$, because $\sigma|_{U_6}=0$. An alternating form on a six-dimensional space with a four-dimensional isotropic subspace has rank at most four. Thus $s_{U_8}|_{\PP(A_4)}=0$, which gives the scheme-theoretic inclusion $\PP(A_4)\subset Q_{U_8}$. At the generic point of $\Gamma$, write $U_7$ for $(\mathcal V_7)_K\otimes_K\CC(\Gamma)$. On $D\times_S\operatorname{Spec}\CC(\Gamma)$, let $\mathcal U_6$ be the pullback of the tautological bundle on $\operatorname{Gr}(6,V_{10})$. The quotient by $U_7$ induces a map
	$$
		\psi:\mathcal U_6\longrightarrow
		\bigl((V_{10}\otimes_{\CC}\CC(\Gamma))/U_7\bigr)
		\otimes_{\CC(\Gamma)}\cO_{D\times_S\operatorname{Spec}\CC(\Gamma)}.
	$$
	The map $\psi$ has rank one at every geometric point. Hence
	$$
		\mathcal A_4:=\ker(\psi)/(V_1\otimes\cO)
	$$
	is a rank-four subbundle of $(U_7/V_1)\otimes\cO$, with fiber $(U_6\cap U_7)/V_1$ at $[U_6]$. The inclusion $\PP(A_4)\subset Q_{U_8}$ on each fiber shows that $\mathcal A_4$ defines the morphism~\eqref{pes:conic-to-fano}.
\end{proof}

\begin{lemma}\label{pes:singular-members}
	In the setup of Section~\ref{pes:geometric-setup}, let $K=\CC(B)$ and let $\Pi_K=\PP_K(\mathcal N_K)$ be the pencil defining $X_{1.2/K}^\sigma$. Under the identification of Lemma~\ref{pes:incidence-variety},
	$$
		\Gamma_K=\Gamma\times_B\operatorname{Spec} K\subset\Pi_K
	$$
	is the reduced closed subscheme parametrizing the singular members of the pencil. Over the algebraic closure $\overline K$, there are exactly three singular members, each of corank two.
\end{lemma}

\begin{proof}
	By Lemma~\ref{pes:incidence-variety}, write
	$$
		\Gamma_K\times_K\overline K=\{p_1,p_2,p_3\}
		\subset\Pi_{\overline K}
	$$
	with $p_1,p_2,p_3$ distinct. For each $i$, set $s_i=\pi_S(p_i)$ and let $q_{p_i}$ be a quadratic form defining the pencil member indexed by $p_i$. Since $D_{s_i}$ is a conic over $\overline K$, it has a $\overline K$-point. Lemma~\ref{pes:conic-spaces} gives a four-dimensional totally isotropic subspace
	$$
		A_{4,i}\subset(\mathcal V_7/V_1)_K\otimes_K\overline K
	$$
	for the quadratic form $q_{p_i}$. Since $A_{4,i}\subset A_{4,i}^{\perp_{q_{p_i}}}$, we have
	$$
		4\leq\dim A_{4,i}^{\perp_{q_{p_i}}}
		=2+\dim\bigl(A_{4,i}\cap\operatorname{rad}(q_{p_i})\bigr).
	$$
	Hence $\operatorname{corank}(q_{p_i})\geq2$ for $i=1,2,3$. We choose a basis of $\mathcal N_K$ and a trivialization of $\det(\mathcal V_7/V_1)_K^\vee$. By Lemma~\ref{pes:pfaffian}, these choices give generators $q_0,q_1$ of the pencil. In a basis of $(\mathcal V_7/V_1)_K$, let $M_0,M_1\in\operatorname{Mat}_6(K)$ be the symmetric matrices of their polar forms. The member with parameter $[s:t]\in\Pi_K\simeq\PP_K^1$ has polar matrix
	$$
		M(s,t)=sM_0+tM_1.
	$$
	Set
	$$
		\Delta(s,t)=\det M(s,t)
		\in H^0(\Pi_K,\cO_{\Pi_K}(6)).
	$$
	The generic member is nondegenerate, so $\Delta\ne0$ and $\deg\operatorname{div}(\Delta)=6$. Fix a closed point $p\in\Pi_{\overline K}$. Let $R_p=\cO_{\Pi_{\overline K},p}$, choose a uniformizer $z_p\in R_p$, and trivialize $\cO_{\Pi_{\overline K}}(1)$ near $p$. The matrix $M(s,t)$ then gives a matrix $M_p\in\operatorname{Mat}_6(R_p)$. Write $M_p(p)=M_p\bmod(z_p)\in\operatorname{Mat}_6(\overline K)$ and put $c_p=6-\operatorname{rank}M_p(p)$. Choose $P_p\in\mathrm{GL}_6(\overline K)$ so that the last $c_p$ rows of $P_pM_p(p)$ vanish. The last $c_p$ rows of $P_pM_p$ are divisible by $z_p$, so
	$$
		P_pM_p=
		\begin{pmatrix}
			I_{6-c_p}&0\\
			0&z_pI_{c_p}
		\end{pmatrix}N_p
		\qquad\text{for some }N_p\in\operatorname{Mat}_6(R_p).
	$$
	Taking determinants gives
	$$
		\operatorname{ord}_p(\Delta)
		=c_p+\operatorname{ord}_{z_p}(\det N_p)\geq c_p.
	$$
	In particular,
	$$
		6=\deg\operatorname{div}(\Delta)
		\geq\sum_{i=1}^3\operatorname{ord}_{p_i}(\Delta)
		\geq\sum_{i=1}^3c_{p_i}\geq6.
	$$
	All inequalities are equalities. Thus $c_{p_i}=\operatorname{ord}_{p_i}(\Delta)=2$ for $i=1,2,3$, and
	$$
		\operatorname{div}(\Delta)=2[p_1]+2[p_2]+2[p_3]
		\quad\text{on }\Pi_{\overline K}.
	$$
	The singular parameters are the zeros of $\Delta$. Hence $V_{\Pi_K}(\Delta)_{\mathrm{red}}$ and $\Gamma_K$ become equal after base change to $\overline K$, and therefore are equal over $K$.
\end{proof}

\subsection{The rationality of \texorpdfstring{$X_1^\sigma$}{X1(sigma)}}
We set
$$
	H=c_1(\mathcal U_1^\vee|_{X_{1.2}^\sigma}),\qquad
	h=\varpi^*c_1(\cO_B(1)).
$$
These first Chern classes act on the Chow groups of $X_{1.2}^\sigma$ even when it is singular. The following is a cycle-theoretic version of \cite[Theorem~3.21]{BF24}, with no geometric integrality assumption on the generic curves. It follows from the argument in Section~\ref{cyc:section} and the projection formula.

\begin{corollary}\label{pes:cycle-criterion}
	In the setup of Section~\ref{pes:geometric-setup}, if a cycle with integer coefficients $z\in\operatorname{CH}_4(X_{1.2}^\sigma)$ satisfies
	$$
		\deg(Hh^3\cap z)\equiv1\pmod2,
	$$
	then $X_{1.2/K}^\sigma$ contains a $K$-line and is $K$-rational. Moreover, $X_1^\sigma$ is rational over $\CC$.
\end{corollary}

\begin{proof}
	Apply Theorem~\ref{cyc:relative} to $f=\varpi:X_{1.2}^\sigma\to B$, with $d=3$, $r=1$, $n=5$, $\mathcal H=\cO_{\PP_B(\mathcal V_7/V_1)}(1)|_{X_{1.2}^\sigma}$, and $\mathcal A_1=\mathcal A_2=\mathcal A_3=\cO_B(1)$. The theorem gives a $K$-line in $X_{1.2/K}^\sigma$ and a birational map $X_{1.2}^\sigma\dashrightarrow B\times\PP^3$. Since $B\simeq\PP^3$ and $X_1^\sigma$ is birational to $X_{1.2}^\sigma$, the variety $X_1^\sigma$ is rational.
\end{proof}

\begin{theorem}\label{pes:rationality}\label{pes:beta-rational}
	Let $\sigma\in U$, where $U$ is the open subset fixed in Assumption~\ref{pes:generic}. If $\beta=0$, then the generic fiber $X_{1.2/K}^\sigma$ contains a $K$-line $L\subset X_{1.2/K}^\sigma$. The closure $Z_L\subset X_{1.2}^\sigma$ of $L$ is a four-dimensional subvariety satisfying
	$$
		\deg(Hh^3\cap[Z_L])=1.
	$$
	In particular, $X_1^\sigma$ is rational.
\end{theorem}

\begin{proof}
	The assumption $\beta=0$ gives a rational section $r:S\dashrightarrow D$ of $\pi$. Since $\pi_S:\Gamma\to S$ is a $\PP^1$-bundle, the generic point of $\Gamma$ maps to the generic point of $S$, where $r$ is defined. Thus $r$ induces a $\CC(\Gamma)$-point of $D\times_S\operatorname{Spec}\CC(\Gamma)$. By Lemma~\ref{pes:incidence-variety} and Lemma~\ref{pes:conic-spaces}, this point determines a four-dimensional subspace $A_4\subset(\mathcal V_7/V_1)_K\otimes_K \CC(\Gamma)$ whose projectivization lies in the member of the pencil indexed by $\Gamma_K=\operatorname{Spec}\CC(\Gamma)$. Let $\iota_1,\iota_2,\iota_3:\CC(\Gamma)\hookrightarrow K^{\mathrm{sep}}$ be the three $K$-embeddings. Write $A_{4,i}=A_4\otimes_{\CC(\Gamma),\iota_i}K^{\mathrm{sep}}$ and let $Q_i$ be the corresponding member of the pencil. Then $\PP(A_{4,i})\subset Q_i$, and the group $\operatorname{Gal}(K^{\mathrm{sep}}/K)$ permutes the pairs $(Q_i,A_{4,i})$. We apply the complementary-pair argument from the proof of Theorem~\ref{cub:cubic-descent}. The subgroup $\operatorname{Gal}(K^{\mathrm{sep}}/\iota_1(\CC(\Gamma)))$ permutes $A_{4,2}$ and $A_{4,3}$, so $A_{4,2}\cap A_{4,3}$ descends to a $\CC(\Gamma)$-subspace $A\subset(\mathcal V_7/V_1)_K\otimes_K \CC(\Gamma)$. Moreover,
	$$
		\dim_{\CC(\Gamma)}A=\dim(A_{4,2}\cap A_{4,3})\geq4+4-6=2,
	$$
	and~\eqref{pes:base-scheme} gives
	$$
		\PP(A_{4,2}\cap A_{4,3})\subset Q_2\cap Q_3=X_{1.2/K^{\mathrm{sep}}}^\sigma.
	$$
	Hence $\PP_{\CC(\Gamma)}(A)\subset X_{1.2/\CC(\Gamma)}^\sigma$. Choose a two-dimensional $\CC(\Gamma)$-subspace $A_2\subset A$ and set $L_{\CC(\Gamma)}=\PP_{\CC(\Gamma)}(A_2)$. We choose generators $q_0,q_1$ of the pencil defining $X_{1.2/K}^\sigma$. Since $L_{\CC(\Gamma)}\subset X_{1.2/\CC(\Gamma)}^\sigma$, the space $A_2\otimes_{\CC(\Gamma)}\CC(\Gamma)(t)$ is totally isotropic for $(q_0+tq_1)_{\CC(\Gamma)(t)}$. Here, $(q_0+tq_1)_{\CC(\Gamma)(t)}$ is the scalar extension of $q_0+tq_1$ to the function field $\mathbf{C}(\Gamma)(t)$. Therefore
	$$
		i_W\bigl((q_0+tq_1)_{\CC(\Gamma)(t)}\bigr)\geq2.
	$$
	The extension $\CC(\Gamma)(t)/K(t)$ has degree three. Theorem~\ref{thm.Springer}, applied to the anisotropic part of the nondegenerate form $q_0+tq_1$, gives
	$$
		i_W(q_0+tq_1)
		=i_W\bigl((q_0+tq_1)_{\CC(\Gamma)(t)}\bigr)\geq2.
	$$
	Theorem~\ref{thm.Amer} now gives a $K$-line $L\subset X_{1.2/K}^\sigma$. One can also obtain $L$ from Proposition~\ref{cyc:linear-space}: for the finite projection $\nu:X_{1.2/\CC(\Gamma)}^\sigma\to X_{1.2/K}^\sigma$, the one-cycle $\nu_*[L_{\CC(\Gamma)}]$ has degree $[\CC(\Gamma):K]=3$. Let $Z_L\subset X_{1.2}^\sigma$ be the closure of $L$. The variety $Z_L$ is integral of dimension $\dim B+\dim L=4$. Since $\deg_K\cO_L(1)=1$, proper pushforward gives
	$$
		\varpi_*(H\cap[Z_L])=[B].
	$$
	The projection formula yields
	$$
		\begin{aligned}
			\deg(Hh^3\cap[Z_L])
			&=\deg\bigl(c_1(\cO_B(1))^3\cap\varpi_*(H\cap[Z_L])\bigr)\\
			&=\deg\bigl(c_1(\cO_B(1))^3\cap[B]\bigr)=1.
		\end{aligned}
	$$
	Since the geometric nonsmooth locus of $X_{1.2/K}^\sigma$ is finite, Proposition~\ref{cyc:projection} gives a birational map $X_{1.2/K}^\sigma\dashrightarrow\PP_K^3$ by projection from $L$. It follows that $X_{1.2}^\sigma$ is birational over $B$ to $B\times\PP^3$. Since $B\simeq\PP^3$ and $X_1^\sigma$ is birational to $X_{1.2}^\sigma$, the variety $X_1^\sigma$ is rational.
\end{proof}

\subsection{Nonemptiness of the locus \texorpdfstring{$\beta=0$}{beta=0}}\label{pes:nonvacuous-section}

The hypothesis $\beta=0$ in Theorem~\ref{pes:beta-rational} holds at some point of the prescribed open subset $U$. We deduce this from the lattice computations of \cite[\S5.2]{BS21}, the criterion of \cite[Proposition~3.3]{BF24}, and the period argument in Appendix~\ref{app:nonvacuous}.

\begin{proposition}
	\label{pes:nonvacuous}
	Let $U\subset\mathcal D^{1,6,10}$ be the nonempty Zariski open subset fixed in Assumption~\ref{pes:generic}. There exists $\sigma\in U$ for which the Brauer class $\beta$ of the associated conic fibration $\pi:D\to S$ vanishes.
\end{proposition}

\begin{proof}
	Lemma~\ref{pes:odd-class} gives a point $\sigma\in U'\subset U$ and a class $C\in\operatorname{NS}(X_6^\sigma)$ with $q(C,[D])=1$. Choose a line bundle $\mathcal L$ with $c_1(\mathcal L)=C$. By Lemma~\ref{pes:odd-pairing}, the restriction $i^*\mathcal L$ to $D$ has degree one on each fiber of $\pi$, where $i:D\hookrightarrow X_6^\sigma$ is the inclusion. Thus \cite[Proposition~3.3]{BF24} gives $\beta=0$.
\end{proof}

\section{Clifford categories and a nodal curve}
\label{der:section}

Throughout this section, $K$ is a field of characteristic zero.  We first prove a decomposition for  pencils of quadrics and then apply it to the generic fiber of $\varpi:X_{1.2}^{\sigma}\to B$.

\subsection{The Clifford component without smoothness assumptions}

Let $E$ be a six-dimensional $K$-vector space, and let $q_0,q_1:E\to K$ be linearly independent quadratic forms such that
$$
	X := V(q_0,q_1)\subset\PP_K(E)
$$
is a complete intersection.  Assume that $q_0+tq_1$ is nondegenerate over $K(t)$. Since $K$ is infinite, we may replace the basis of the pencil so that both $q_0$ and $q_1$ are nondegenerate. Set $\Pi=\PP_K(\langle q_0,q_1\rangle)$, and let
$$
	\rho:\mathscr Q\longrightarrow\Pi
$$
be the universal quadric of the pencil.  Thus $\mathscr Q$ is the hypersurface in $\PP_K(E)\times_K\Pi$ with equation $aq_0+bq_1=0$. Write $\mathscr C_0$ for the even Clifford algebra of $\rho$. The morphism $\rho$ is flat: every member is a nonzero quadratic form, and hence every fiber is a degree-two hypersurface in $\PP_K(E)$.

Homological projective duality gives a fully faithful functor
\begin{equation}\label{der:clifford-functor}
	\Phi:\Db(\Pi,\mathscr C_0)\longrightarrow\Db(X)
\end{equation}
and a semiorthogonal decomposition
\begin{equation}\label{der:clifford-sod}
	\Db(X)=\left\langle
	\Phi\Db(\Pi,\mathscr C_0),\cO_X,\cO_X(1)
	\right\rangle.
\end{equation}
Here we use \cite[Theorem~2.19(1)]{ABB14} over $\Spec K$, and tensor its decomposition by $\cO_X(-1)$.  This is the version over $K$ of \cite[Theorem~5.5]{Kuz08Q}.  The hypothesis on the linear section is that $X$ is a complete intersection; smoothness of $X$ is not required. The exceptional pair in \eqref{der:clifford-sod} remains exceptional for a singular complete intersection: the Koszul resolution on $\PP_K(E)$ gives
$$
	\mathbf{R}\Gamma(X,\cO_X)=K,
	\qquad \mathbf{R}\Gamma(X,\cO_X(-1))=0.
$$

\begin{proposition}\label{der:relative-clifford}
	Let $U$ be a smooth variety over $K$.  Let $\mathcal E$ be a rank-six vector bundle on $U$, and let $\mathcal L\subset\operatorname{Sym}^2\mathcal E^\vee$ be a rank-two subbundle which is the direct sum of two line subbundles defining quadratic forms that are nondegenerate over the function field of $U$. Suppose that the common zero scheme
	$$
		f_U:X_U\longrightarrow U,
		\qquad X_U\subset\PP_U(\mathcal E),
	$$
	is a flat relative complete intersection of two quadrics.  Let $\pi_U:\Pi_U=\PP_U(\mathcal L)\to U$ be the parameter space of the pencil, and let $\rho_U:\mathscr Q_U\to\Pi_U$ be its universal quadric, with even Clifford algebra $\mathscr C_{0,U}$.  Then there is a fully faithful functor
	$$
		\Phi_U:\Db(\Pi_U,\mathscr C_{0,U})\longrightarrow\Db(X_U)
	$$
	and a semiorthogonal decomposition
	\begin{equation}\label{der:relative-sod}
		\Db(X_U)=\left\langle
		\Phi_U\Db(\Pi_U,\mathscr C_{0,U}),
		Lf_U^*\Db(U),
		Lf_U^*\Db(U)\otimes\cO_{X_U}(1)
		\right\rangle.
	\end{equation}
	In particular, this conclusion applies on a sufficiently small dense open subset of the base $B$ in the Peskine construction.
\end{proposition}

\begin{proof}
	Apply relative homological projective duality in \cite[Theorem~2.19(1)]{ABB14} to the rank-six bundle $\mathcal E$ and the rank-two pencil $\mathcal L$. Tensoring the resulting decomposition by $\cO_{X_U}(-1)$ gives~\eqref{der:relative-sod}. For $j\in\{0,1\}$, let
	$$
		F_j:\Db(U)\longrightarrow\Db(X_U),
		\qquad G\longmapsto Lf_U^*G\otimes\cO_{X_U}(j).
	$$
	Flatness of $f_U$ makes both functors $F_j$ well defined. The relative Koszul resolution gives
	$$
		Rf_{U*}\cO_{X_U}=\cO_U,
		\qquad Rf_{U*}\cO_{X_U}(-1)=0.
	$$
	The first identity and the projection formula show that $F_0$ and $F_1$ are fully faithful. The second gives $\RHom_{X_U}(F_1(G),F_0(H))=0$ for all $G,H\in\Db(U)$. For the last assertion, apply the proposition over $\CC$. Choose nondegenerate generators $q_0,q_1$ of the Peskine pencil on $(\mathcal V_7/V_1)_K$, where $K=\CC(B)$. Extend $(\mathcal V_7/V_1)_K$ to a rank-six bundle $\mathcal E$ over a dense open subset $U\subset B$, and extend $q_0,q_1$ to sections of $\Sym^2\mathcal E^\vee$. After shrinking $U$, the sections $q_0,q_1$ generate a rank-two subbundle
	$$
		\mathcal L=\cO_Uq_0\oplus\cO_Uq_1
		\subset\Sym^2\mathcal E^\vee,
	$$
	and $X_U=V(q_0,q_1)\subset\PP_U(\mathcal E)$ is flat over $U$, with every fiber of codimension two. Since $U\subset B=\PP^3_{\CC}$ is smooth, the hypotheses of the proposition hold. No smoothness of $X_U$ or of $f_U$ is used.
\end{proof}

\subsection{A criterion for passing to bounded derived categories}

The geometric description of a Clifford component that we use below is initially an equivalence of unbounded categories.  We give the boundedness argument separately, since an argument for perfect complexes alone does not give an equivalence of bounded derived categories of coherent sheaves when the varieties are singular.

\begin{lemma}\label{der:bounded-test}
	Let $Z$ be a projective $K$-scheme and let $\mathscr A$ be an Azumaya algebra on $Z$.  For $F\in\Dqc(Z,\mathscr A)$, the following conditions are equivalent:
	\begin{enumerate}
		\item $F\in\Db(Z,\mathscr A)$;
		\item for every $P\in\operatorname{Perf}(Z,\mathscr A)$, the complex
		$$
			\ \mathbf{R} \mathrm{Hom}_{(Z,\mathscr A)}(P,F)
		$$
		belongs to $D^b(\mathrm{Vect}(K))$, the bounded derived category of finite-dimensional $K$-vector spaces.
	\end{enumerate}
	The statement also holds with $\mathscr A=\cO_Z$.
\end{lemma}

\begin{proof}
	Suppose that $F\in\Db(Z,\mathscr A)$ and $P\in\Perf(Z,\mathscr A)$. Then $\mathbf R\mathcal{H}om_{\mathscr A}(P,F)$ is a bounded coherent complex on $Z$. Since $Z$ is proper,
	$$
		\RHom_{(Z,\mathscr A)}(P,F)
		\simeq\mathbf R\Gamma\bigl(Z,\mathbf R\mathcal{H}om_{\mathscr A}(P,F)\bigr)
		\in D^b(\mathrm{Vect}(K)).
	$$
    Conversely, choose a closed immersion $i:Z\hookrightarrow\PP_K^N$ and put $\cO_Z(1)=i^*\cO_{\PP_K^N}(1)$. For $0\leq j\leq N$, the right $\mathscr A$-module $\mathscr A\otimes\cO_Z(j)$ is perfect. Write $F_{\cO_Z}$ for the underlying complex of $\cO_Z$-modules. The hypothesis gives
	$$
		\RHom_{(Z,\mathscr A)}
		(\mathscr A\otimes\cO_Z(j),F)
		\simeq \mathbf R\Gamma(Z,F_{\cO_Z}(-j))
		\in D^b(\mathrm{Vect}(K)).
	$$
	Beilinson's resolution of the diagonal \cite{Bei79} constructs $Ri_*F_{\cO_Z}$ by finitely many cones and shifts from the complexes
	$$
		\mathbf R\Gamma(Z,F_{\cO_Z}(-j))
		\otimes_K\Omega^j_{\PP_K^N}(j),
		\qquad 0\leq j\leq N.
	$$
	Each of these complexes is bounded coherent, so $Ri_*F_{\cO_Z}\in\Db(\PP_K^N)$. Since $i_*$ is exact on quasi-coherent sheaves and detects coherence, $F_{\cO_Z}$ is bounded coherent. The algebra $\mathscr A$ is finite locally free over $\cO_Z$, so $F\in\Db(Z,\mathscr A)$.
\end{proof}

\begin{lemma}\label{der:residual-bounded}
	Let $Z$ be a projective $K$-scheme with a finite exceptional collection $E_1,\ldots,E_m$ of perfect complexes.  Let
	$$
		\mathcal K_{\mathrm{qc}}=
		\{F\in\Dqc(Z)\mid
		\RHom_Z(E_i,F)=0\text{ for every }1\leq i\leq m\}.
	$$
	Suppose $M$ is a projective $K$-scheme, $\gamma\in\rBr(M)$, and there is a $K$-linear exact equivalence
	$$
		F_M:\Dqc(M,\gamma)\xrightarrow{\sim}\mathcal K_{\mathrm{qc}}.
	$$
	Then $F_M$ restricts to an equivalence
	$$
		F_M:\Db(M,\gamma)\xrightarrow{\sim}
		\mathcal K_{\mathrm{qc}}\cap\Db(Z).
	$$
\end{lemma}

\begin{proof}
	For $F\in\Dqc(Z)$, set $F_m=F$. For $i=m,m-1,\ldots,1$, define $F_{i-1}$ by the evaluation triangle
	$$
		\RHom_Z(E_i,F_i)\otimes_K E_i
		\longrightarrow F_i\longrightarrow F_{i-1}.
	$$
	Then $F_0$ belongs to the right orthogonal $\mathcal K_{\mathrm{qc}}$, and $F\mapsto F_0$ is the projection
	$$
		\operatorname{pr}_{\mathcal K}:\Dqc(Z)
		\longrightarrow\mathcal K_{\mathrm{qc}}.
	$$
	For $P\in\Perf(Z)$, each complex $\RHom_Z(E_i,P)$ is bounded and finite-dimensional. The evaluation triangles therefore show that $\operatorname{pr}_{\mathcal K}$ preserves perfect complexes. Since every $E_i$ is perfect, $\operatorname{pr}_{\mathcal K}$ also preserves coproducts. The objects $\operatorname{pr}_{\mathcal K}P$, with $P\in\Perf(Z)$, compactly generate $\mathcal K_{\mathrm{qc}}$. Their thick closure is the category of compact objects, so
	$$
		\mathcal K_{\mathrm{qc}}^c
		=\mathcal K_{\mathrm{qc}}\cap\Perf(Z).
	$$
    For $P\in\Perf(Z)$, the evaluation triangles give a triangle
	$$
		A_P\longrightarrow P\longrightarrow\operatorname{pr}_{\mathcal K}P,
		\qquad A_P\in\langle E_1,\ldots,E_m\rangle.
	$$
	If $G\in\mathcal K_{\mathrm{qc}}$, then $\RHom_Z(A_P,G)=0$. Applying $\RHom_Z(-,G)$ to this triangle yields
	$$
		\RHom_Z(P,G)
		\simeq
		\RHom_{\mathcal K_{\mathrm{qc}}}
		(\operatorname{pr}_{\mathcal K}P,G).
	$$
	By Lemma~\ref{der:bounded-test}, $G$ belongs to $\Db(Z)$ if and only if $\RHom_{\mathcal K_{\mathrm{qc}}}(Q,G)\in D^b(\mathrm{Vect}(K))$ for every compact object $Q\in\mathcal K_{\mathrm{qc}}$. The equivalence $F_M$ preserves compact objects and all graded homomorphism spaces. Applying Lemma~\ref{der:bounded-test} on $M$ therefore gives
	$$
		F\in\Db(M,\gamma)
		\quad\Longleftrightarrow\quad
		F_M(F)\in\mathcal K_{\mathrm{qc}}\cap\Db(Z),
	$$
	as required.
\end{proof}

\subsection{Pencils with three corank-two members}

Return to the pencil $\rho:\mathscr Q\to\Pi$ defined above.  Assume that it has a nondegenerate member and that its geometric singular members occur at exactly three distinct parameters, each with corank two. Let
$$
	\Delta_2\hookrightarrow\Pi
$$
denote the corank-two locus.  It is a reduced finite \mbox{\'{e}tale} $K$-scheme of degree three.  For the scheme-theoretic assertion, note that the determinant has degree six and that its vanishing order is at least two at every such parameter.  All three orders are therefore exactly two.  Over a local discrete valuation ring at a geometric parameter, the two nonunit elementary divisors of the matrix have valuation one. Its corank-two determinantal ideal is consequently generated by a uniformizer.

There is a finite \mbox{\'{e}tale} double cover
\begin{equation}\label{der:component-cover}
	\epsilon:\widetilde\Delta_2\longrightarrow\Delta_2
\end{equation}
obtained from the Stein factorization of $\mathcal{M}_3(\mathscr Q/\Pi)|_{\Delta_2}\to\Delta_2$. Its two geometric points distinguish the two families of maximal projective three-spaces in a corank-two quadric \cite[Lemma~4.12]{CO26}.

\begin{lemma}\label{der:cover-splitting}
	Suppose that for every closed point $s\in\Delta_2$, there exist a smooth projective geometrically integral conic $D_s$ over $K(s)$ and a $K(s)$-morphism
	$$
		a_s:D_s\longrightarrow
		\mathcal{M}_3(\mathscr Q_s/K(s)).
	$$
	Then the double cover $\epsilon$ in \eqref{der:component-cover} has a section. No $K(s)$-rational point on $D_s$ is required.
\end{lemma}

\begin{proof}
	Fix a closed point $s\in\Delta_2$. Compose $a_s$ with the morphism from the Stein factorization to obtain
	$$
		D_s\xrightarrow{a_s}\mathcal{M}_3(\mathscr Q_s/K(s))
		\longrightarrow\widetilde\Delta_{2,s}.
	$$
	The scheme $\widetilde\Delta_{2,s}$ is affine, and $H^0(D_s,\cO_{D_s})=K(s)$. Hence this composite factors as
	$$
		D_s\longrightarrow\Spec K(s)
		\xrightarrow{\eta_s}\widetilde\Delta_{2,s}.
	$$
	The morphism $\eta_s$ is a section of $\epsilon$ over $s$. Since
	$$
		\Delta_2=\coprod_{s\in|\Delta_2|}\Spec K(s),
	$$
	the sections $\eta_s$ give a section of $\epsilon:\widetilde\Delta_2\to\Delta_2$.
\end{proof}

We apply \cite[Theorem~7.25]{CO26} to this pencil. The passage to bounded derived categories and the description of the nodal curve are proved below. For smooth complete intersections of two quadrics over an arbitrary field, a decomposition of this type is \cite[8.7(i)]{CO26}; here $X$ is singular. As remarked in \cite{CO26}, statements for bounded derived categories can also be deduced by adapting \cite[Theorem~6.2]{BS20}; Lemma~\ref{der:residual-bounded} gives a direct argument in our setting.

Recall that the dual graph of a nodal curve has one vertex for each irreducible component and one edge for each node, joining the vertices of the components containing its two branches, with loops allowed. For a connected projective nodal curve $C$ over an algebraically closed field, with $n$ irreducible components and $\delta$ nodes, the normalization exact sequence gives $p_a(C)=\sum_{i=1}^n g(\widetilde C_i)+\delta-n+1$, where $\widetilde C_i$ is the normalization of the $i$-th component; for example, see \cite[\S3.1]{Cap18}. The complete bipartite graph $K_{2,3}$ has two sets of vertices, of sizes two and three, with an edge between every pair from different sets and no edges within either set.

\begin{theorem}\label{der:curve-model}
	Let $X\subset\PP_K(E)$ be the complete intersection above.  Assume that its pencil has exactly three geometric singular members, all of corank two, and that the double cover $\epsilon:\widetilde\Delta_2\to\Delta_2$ has a section.  Then there exist a projective nodal $K$-curve $M$, a proper morphism $\mu:M\to\Pi$, a Brauer class $\gamma\in\rBr(M)$, and a fully faithful functor
	$$
		\Psi:\Db(M,\gamma)\longrightarrow\Db(X)
	$$
	such that
	\begin{equation}\label{der:curve-sod}
		\Db(X)=\left\langle
		\Psi\Db(M,\gamma),\cO_X,\cO_X(1)
		\right\rangle.
	\end{equation}
	The morphism $\mu$ is finite of degree two away from $\Delta_2$ and has a geometric fiber $\PP^1$ over each point of $\Delta_2$. Over $\overline K$, the curve $M$ has five irreducible components, each isomorphic to $\PP^1_{\overline K}$, and six nodes.  Its dual graph is the complete bipartite graph $K_{2,3}$, and $p_a(M_{\overline K})=2$.
\end{theorem}

\begin{proof}
	Define the unbounded residual category of the universal quadric $\rho$ by
	$$
		\mathcal K_{\mathrm{qc}}(\mathscr Q/\Pi)
		=\{F\in\Dqc(\mathscr Q)\mid
		R\rho_*(F\otimes\cO_{\mathscr Q}(-j))=0
		\text{ for }j=0,1,2,3\}.
	$$
	The quadric fibration $\rho$ has relative dimension four and corank at most two. Its generic fiber is nondegenerate, the base $\Pi$ is integral, and the double cover~\eqref{der:component-cover} is split. Thus \cite[Theorem~7.25]{CO26} gives an algebraic space $M$, a twist $\gamma$, and a $\Pi$-linear equivalence
	\begin{equation}\label{der:unbounded-equivalence}
		\Dqc(M,\gamma)\xrightarrow{\sim}
		\mathcal K_{\mathrm{qc}}(\mathscr Q/\Pi).
	\end{equation}
	The theorem also gives a proper morphism $\mu:M\to\Pi$ that is finite of degree two over $\Pi\setminus\Delta_2$ and has geometric fiber $\PP^1$ over each point of $\Delta_2$. Hence $M$ is a proper one-dimensional algebraic space over $K$. A field is a complete Noetherian local ring, so \cite[Tag~0AE7]{Stacks} shows that $M$ is a projective $K$-scheme. Let $\mathfrak G\to M$ be the gerbe in the construction of~\eqref{der:unbounded-equivalence}, and let $\mathscr S$ be the universal dual spinor sheaf on $\mathscr Q\times_\Pi\mathfrak G$. Write $p:\mathscr Q\times_\Pi\mathfrak G\to\mathfrak G$ for the projection. The sheaf $\mathscr S$ is flat over $\mathfrak G$ by \cite[\S7.1]{CO26}. For a sufficiently large integer $N$, relative Serre vanishing and cohomology and base change show that
	$$
		\mathcal V=p_*\bigl(\mathscr S\otimes\cO_{\mathscr Q}(N)\bigr)
	$$
	is locally free of positive rank; here $\cO_{\mathscr Q}(N)$ is pulled back to $\mathscr Q\times_\Pi\mathfrak G$. The endomorphism algebra $\mathcal{E}nd(\mathcal V)$ descends to an Azumaya algebra on $M$. Thus the twist in~\eqref{der:unbounded-equivalence} is a Brauer class. We write $\gamma$ for this class, with its sign fixed by the right-module convention in~\eqref{der:unbounded-equivalence}.

	The quadric-fibration decomposition has four further components
	$$
		\mathbf L\rho^*\Dqc(\Pi)\otimes\cO_{\mathscr Q}(j),
		\qquad j=0,1,2,3.
	$$
	The full exceptional pair $\langle\cO_\Pi,\cO_\Pi(1)\rangle$ gives the perfect exceptional objects
	$$
		E_{j,e}=\mathbf L\rho^*\cO_\Pi(e)\otimes\cO_{\mathscr Q}(j),
		\qquad j=0,1,2,3,\quad e=0,1.
	$$
	Order the objects $E_{j,e}$ first by increasing $j$, and then by increasing $e$. Since $\cO_\Pi$ and $\cO_\Pi(1)$ compactly generate $\Dqc(\Pi)$, the right orthogonal of these eight objects is $\mathcal K_{\mathrm{qc}}(\mathscr Q/\Pi)$. The scheme $\mathscr Q$ is projective over $K$, so Lemma~\ref{der:residual-bounded} restricts~\eqref{der:unbounded-equivalence} to
	\begin{equation}\label{der:bounded-equivalence}
		\Db(M,\gamma)\xrightarrow{\sim}
		\mathcal K_{\mathrm{qc}}(\mathscr Q/\Pi)\cap\Db(\mathscr Q).
	\end{equation}
	The quadric-fibration decomposition \cite[Theorem~2.11]{ABB14} extends \cite[Theorem~4.2]{Kuz08Q} to the base field $K$. Tensoring that decomposition by $\cO_{\mathscr Q}(-1)$ gives
	$$
		\mathcal K_{\mathrm{qc}}(\mathscr Q/\Pi)\cap\Db(\mathscr Q)
		\simeq\Db(\Pi,\mathscr C_0).
	$$
	Composing this equivalence with~\eqref{der:bounded-equivalence} and $\Phi$ defines
	$$
		\Psi:\Db(M,\gamma)\xrightarrow{\sim}\Db(\Pi,\mathscr C_0)
		\xrightarrow{\Phi}\Db(X).
	$$
	The functor $\Psi$ is fully faithful, and~\eqref{der:clifford-sod} gives~\eqref{der:curve-sod}. To describe $M_{\overline K}$, choose a nondegenerate member $q_0$. Let $b_{q_i}$ be the polar form of $q_i$ on $E_{\overline K}$, and define $A\in\End(E_{\overline K})$ by
	$$
		b_{q_1}(v,w)=b_{q_0}(Av,w).
	$$
	Let $s_i=[a_i:b_i]$, for $i=1,2,3$, be the singular parameters. Since $q_0$ is nondegenerate, $b_i\ne0$. Put $\lambda_i=-a_i/b_i$. The corank-two condition gives
	$$
		E_i:=\ker(A-\lambda_i\operatorname{id}),
		\qquad \dim E_i=2,
		\qquad E_{\overline K}=E_1\oplus E_2\oplus E_3.
	$$
	The endomorphism $A$ is self-adjoint for $b_{q_0}$, so the eigenspaces $E_i$ are mutually orthogonal and nondegenerate. Choose coordinates on each $E_i$ such that $q_0|_{E_i}=x_{2i-2}x_{2i-1}$. Set $\ell_i(a,b)=a+b\lambda_i$. Then
	$$
		aq_0+bq_1
		=\ell_1(a,b)x_0x_1\;\perp\;
		\ell_2(a,b)x_2x_3\;\perp\;
		\ell_3(a,b)x_4x_5,
		\qquad V(\ell_i)=\{s_i\}.
	$$ 
    Fix $i\in\{1,2,3\}$. In an \mbox{\'{e}tale} neighborhood of $s_i$, trivialize $\cO_\Pi(1)$ and choose a local parameter $t$ at $s_i$. The coefficient $\ell_i$ is a unit times $t$, and the other two coefficients are units. Rescaling one coordinate in each binary summand therefore gives
	$$
		aq_0+bq_1\simeq\mathbb H\perp\mathbb H\perp t\mathbb H,
	$$
	where $\mathbb H$ is a hyperbolic plane. Removing the first two hyperbolic planes gives the maximal hyperbolic reduction
	\begin{equation}\label{der:local-curve}
		\{tuv=0\}\subset
		\mathbb A^1_t\times\PP^1_{[u:v]}.
	\end{equation}
	By \cite[7.24]{CO26}, $M_{\overline K}$ is \mbox{\'{e}tale} locally isomorphic to~\eqref{der:local-curve}. The three components of this model are the vertical line $V(t)$ and the horizontal branches $V(u)$ and $V(v)$. The line $V(t)$ meets these branches at $(0,[0:1])$ and $(0,[1:0])$. The local equations at the two intersections are $t(u/v)=0$ and $t(v/u)=0$, respectively. Thus both intersections are ordinary nodes, and the model is reduced. Set $\Pi^\circ=\Pi_{\overline K}\setminus\{s_1,s_2,s_3\}$. Let $B_0,B_1$ be the matrices of $b_{q_0},b_{q_1}$ in the chosen coordinates. For some $c\in\overline K^\times$,
	$$
		\det(aB_0+bB_1)=c(\ell_1\ell_2\ell_3)^2.
	$$
	The discriminant cover over $\Pi^\circ$ therefore splits, so
	$$
		\mu^{-1}(\Pi^\circ)\simeq\Pi^\circ\sqcup\Pi^\circ.
	$$
	Let $C_+,C_-$ be the closures of the two sheets in $M_{\overline K}$. The local model~\eqref{der:local-curve} extends each sheet across every $s_i$ as a smooth branch mapped isomorphically to the base. Hence $\mu|_{C_\pm}$ identifies $C_\pm$ with $\Pi_{\overline K}\simeq\PP^1_{\overline K}$. For each $i$, the vertical line in~\eqref{der:local-curve} gives a component $R_i\simeq\PP^1_{\overline K}$ above $s_i$. The intersections are
	$$
		\#(C_+\cap R_i)=\#(C_-\cap R_i)=1,
		\qquad C_+\cap C_-=\varnothing,
		\qquad R_i\cap R_j=\varnothing\quad(i\ne j).
	$$
	All six intersection points are nodes, and there are no other singular points. Thus $M_{\overline K}$ has five rational components and dual graph $K_{2,3}$. The normalization exact sequence gives
	$$
		p_a(M_{\overline K})
		=6-5+1=2.
	$$
\end{proof}

\begin{figure}[ht]
	\centering
	\begin{tikzpicture}[every node/.style={circle,draw,inner sep=3pt}]
		\node (p) at (-1.4,1.1) {$C_+$}; \node (m) at (1.4,1.1) {$C_-$}; \node (r1) at (-2.2,-0.5) {$R_1$}; \node (r2) at (0,-0.5) {$R_2$}; \node (r3) at (2.2,-0.5) {$R_3$}; \draw (p)--(r1) (p)--(r2) (p)--(r3); \draw (m)--(r1) (m)--(r2) (m)--(r3);
	\end{tikzpicture}
	\caption{The dual graph of $M_{\overline K}$. Each vertex represents a projective line, and each edge represents one node of the curve.}
	\label{fig:dual-graph}
\end{figure}
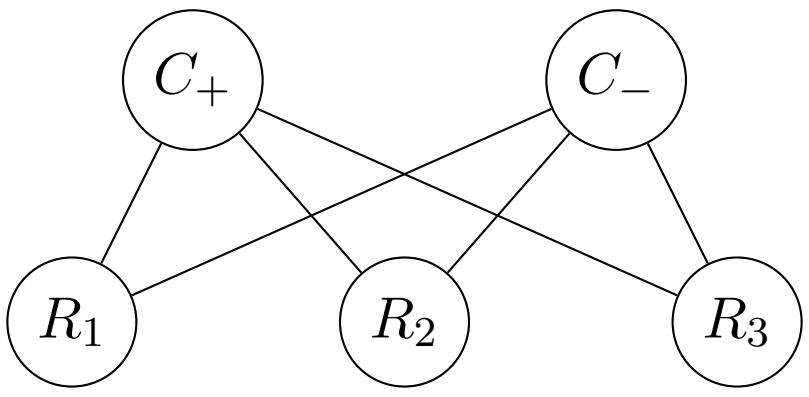

\begin{corollary}\label{der:peskine}
	Let $\sigma$ be as in Section~\ref{pes:geometric-setup}, and put $K=\CC(B)$. Let $X_{1.2/K}^\sigma$ be the generic fiber of $\varpi:X_{1.2}^{\sigma}\to B$.  There exist a projective nodal $K$-curve $M_K$, a Brauer class $\gamma_K\in\rBr(M_K)$, and a fully faithful functor
	$$
		\Psi_K:\Db(M_K,\gamma_K)\longrightarrow\Db(X_{1.2/K}^\sigma)
	$$
	such that
	$$
		\Db(X_{1.2/K}^\sigma)=\left\langle
		\Psi_K\Db(M_K,\gamma_K),\cO_{X_{1.2/K}^\sigma},\cO_{X_{1.2/K}^\sigma}(1)
		\right\rangle.
	$$
	The curve $M_K\times_K\overline K$ has five components, each isomorphic to $\PP^1_{\overline K}$, and six nodes, with dual graph $K_{2,3}$.
\end{corollary}

\begin{proof}
	Let $\rho:\mathscr Q\to\Pi_K$ be the universal quadric of the pencil defining $X_{1.2/K}^\sigma$. By Lemma~\ref{pes:singular-members}, its singular parameters form the degree-three scheme $\Gamma_K$, and every geometric singular member has corank two. Thus $\Delta_2=\Gamma_K$. Fix a closed point $s\in\Gamma_K$, and use $\pi_S:\Gamma\to S$ to define $D_s=D\times_S\Spec K(s)$. Lemma~\ref{pes:conic-spaces} gives the morphism
	$$
		a_s:D_s\longrightarrow\mathcal{M}_3(\mathscr Q_s/K(s)),
		\qquad
		[U_6]\longmapsto[\PP\bigl((U_6\cap U_7)/V_1\bigr)].
	$$
	Here $U_7=(\mathcal V_7)_K\otimes_K K(s)$, and we view $V_1$ as a subspace of $U_7$ after scalar extension. The conic $D_s$ is smooth and geometrically integral, so Lemma~\ref{der:cover-splitting} gives a section of $\epsilon:\widetilde\Delta_2\to\Delta_2$. Theorem~\ref{der:curve-model} now gives the stated decomposition and curve $M_K$. The argument uses the morphism $a_s$, without requiring a $K(s)$-point on $D_s$, so no assumption on $\beta$ is needed.
\end{proof}

\begin{remark}\label{der:scope}
	The class $\gamma_K$ lies in the Brauer group of the curve $M_K$, whereas $\beta$ lies in the Brauer group of the K3 surface $S$. We do not know whether there is a relation between $\gamma_K$ and $\beta$.
\end{remark}

\section{Open questions}
For a general trivector $\sigma$, let $\mathcal A_\sigma$ be the Kuznetsov component of the hyperplane section
$$
	X_3^\sigma = \{[U_3]\in\operatorname{Gr}(3,V_{10})\mid \sigma|_{U_3}=0\}.
$$
Bernardara--Fatighenti--Manivel conjecture that $\Db(X_1^\sigma)$ admits a semiorthogonal decomposition with three components equivalent to $\mathcal A_\sigma$ and four exceptional objects \cite[Conjecture~23]{BFM21}. The semiorthogonal decompositions obtained from even Clifford algebras in the previous section suggest a possible approach to this conjecture through the birational models associated with $\sigma\in\mathcal D^{1,6,10}$.

The first question is whether the decomposition in~\eqref{der:relative-sod} extends from $X_U$ to $X_{1.2}^\sigma$. Such an extension requires an analysis of the fibers over $B\setminus U$. Can the birational correspondence $X_{1.2}^\sigma\dashrightarrow X_1^\sigma$ then be used to construct a semiorthogonal decomposition of $\Db(X_1^\sigma)$ or of a categorical resolution of $X_1^\sigma$? In particular, can this construction produce a component equivalent to $\Db(S,\beta)$? The description of $\Db(X_{1.2/K}^\sigma)$ in Corollary~\ref{der:peskine} does not yet provide a comparison between $\Db(M_K,\gamma_K)$ and $\Db(S,\beta)$.

A further question is whether such a construction can be made compatible with deformations of $\sigma$ out of $\mathcal D^{1,6,10}$, so as to obtain the three components predicted by \cite[Conjecture~23]{BFM21}. The varieties $X_1^\sigma$ and $X_{1.2}^\sigma$ considered here are singular. Thus a decomposition on these special fibers does not directly give a decomposition for general $\sigma$. One would also need to identify the deformed components with $\mathcal A_\sigma$.

Finally, following \cite[\S4]{BF24}, one may ask whether the rationality of $X_1^\sigma$ can be characterized by $\mathcal A_\sigma$. When $X_1^\sigma$ and $X_3^\sigma$ are smooth, is the rationality of $X_1^\sigma$ equivalent to the existence of a smooth projective K3 surface $T$ with
$$
	\mathcal A_\sigma\simeq\Db(T)?
$$
This is analogous to the categorical rationality conjecture for cubic fourfolds \cite[Conjecture~1.1]{Kuz10}. On $\mathcal D^{1,6,10}$, a comparison with $\Db(S,\beta)$ could give a categorical interpretation of Theorem~\ref{pes:beta-rational}. The theorem establishes a sufficient condition for rationality; it does not establish the proposed categorical characterization.

\appendix
\section{Hyperbolic reduction and explicit birational maps}\label{app:hyperbolic}
We use the hypotheses and notation of Proposition~\ref{higher:reduction}. Lemma~\ref{higher:quadric-family} and Lemma~\ref{higher:birational} are essentially the special case of \cite[Lemmas~2.4 and 3.4]{KS18} (see also the proof of \cite[Lemma~3.7]{KS18}) in which the projective dimension of the linear system exceeds the dimension of the linear subspace by one; we include proofs to fix notation.

We choose a decomposition
$$
	K^{n+1}=U\oplus W, \qquad \dim U=m-1, \qquad \dim W=n-m+2.
$$
Fix coordinates $u_1,\ldots,u_{m-1}$ on $U$. Since $q_i|_U=0$, we can write
$$
	q_i(u,w) = \sum_{a=1}^{m-1}u_aL_{ia}(w)+Q_i(w), \qquad 1\leq i\leq m,
$$
where $L_{ia}:W\longrightarrow K$ is a linear form for $1\leq i\leq m$ and $1\leq a\leq m-1$, and $Q_i:W\longrightarrow K$ is quadratic for $1\leq i\leq m$.

For each $w\in W$, let
\begin{equation}\label{eq:M(w)}
	M(w):K^{m-1}\longrightarrow K^m
\end{equation}
be the linear map with matrix $(L_{ia}(w))_{i,a}$, and put
$$
	Q(w)= \begin{pmatrix}
		      Q_1(w)\\
		      \vdots\\
		      Q_m(w)
	       \end{pmatrix}.
$$
Using the chosen coordinates to identify $U$ with $K^{m-1}$, the equations of $Y$ become
$$
	M(w)u+Q(w)=0.
$$
For every nonzero $\lambda=(\lambda_1,\ldots,\lambda_m)\in K^m$, write
$$
	q_\lambda=\sum_{i=1}^m\lambda_iq_i:
	K^{n+1}\longrightarrow K.
$$
The forms $q_1,\ldots,q_m$ are linearly independent because $Y$ has codimension $m$. Thus $[\lambda]\in\PP_K^{m-1}$ parametrizes the members of their linear system. We use the polar form
$$
	b_{q_\lambda}:K^{n+1}\times K^{n+1}\longrightarrow K, \qquad
	b_{q_\lambda}(x,y) = q_\lambda(x+y)-q_\lambda(x)-q_\lambda(y).
$$
In particular,
$$
	b_{q_\lambda}(u,w) = \lambda^{\mathsf t}M(w)u \qquad (u\in U,\ w\in W).
$$
Let
$$
	\operatorname{rad}(q_\lambda) = \{ v \in K^{n+1} \mid b_{q_\lambda}(v,x)=0 \text{ for every } x \in K^{n+1} \}
$$
denote the radical of $q_\lambda$ (Section~\ref{pre:quadratic-forms}). Since the characteristic of $K$ is not $2$, the singular locus of $V(q_\lambda)$ satisfies
$$
	\operatorname{Sing}(V(q_\lambda)) = \mathbf{P}(\operatorname{rad}(q_\lambda)).
$$
All orthogonal complements below are taken with respect to $b_{q_\lambda}$. When passing to an extension field, we retain the same notation for the scalar extensions of the vector spaces and quadratic forms.

\begin{lemma}\label{higher:rank}
	For every nonzero $u\in U$, the linear map
	$$
		W\longrightarrow K^m, \qquad w\longmapsto M(w)u
	$$
	is surjective. Also, for every nonzero $\lambda\in K^m$, the linear map
	$$
		W\longrightarrow U^\vee, \qquad w\longmapsto
		\bigl(u\longmapsto\lambda^{\mathsf t}M(w)u\bigr)
	$$
	is surjective. Moreover,
	$$
		\operatorname{corank}(q_\lambda)\leq m-1.
	$$
	All three assertions remain valid after any extension of $K$.
\end{lemma}

\begin{proof}
	Fix $u\in U\setminus\{0\}$, and define
	$$
		A_u:W\longrightarrow K^m,
		\qquad w\longmapsto M(w)u.
	$$
	At $(u,0)$, the differentials of $q_1,\ldots,q_m$ vanish on $U$ and satisfy
	$$
		\bigl(dq_1|_{(u,0)}(0,w),\ldots,dq_m|_{(u,0)}(0,w)\bigr)^{\mathsf t}
		=A_u(w).
	$$
	Since $Y$ is a smooth complete intersection of codimension $m$, the Jacobian criterion gives $\operatorname{rank}A_u=m$. For $\lambda\in K^m\setminus\{0\}$, define
	$$
		B_\lambda:W\longrightarrow U^\vee,
		\qquad B_\lambda(w)(u)=\lambda^{\mathsf t}M(w)u.
	$$
	If $B_\lambda$ were not surjective, some $u\in U\setminus\{0\}$ would annihilate $\operatorname{Im}B_\lambda$. Then
	$$
		\lambda^{\mathsf t}A_u(w)=0
		\qquad\text{for every }w\in W,
	$$
	contrary to the surjectivity of $A_u$ and the assumption $\lambda\ne0$. For the corank bound, extend scalars to an algebraic closure $\bar K$. The form $q_\lambda$ has the same corank after this extension, and $Y_{\bar K}$ remains smooth. Suppose that $\operatorname{corank}(q_\lambda)\ge m$. After reordering $q_1,\ldots,q_m$, assume that $\lambda_1\ne0$. The identity
	$$
		q_1=\frac{1}{\lambda_1}
		\left(q_\lambda-\sum_{i=2}^m\lambda_iq_i\right)
	$$
	gives
	$$
		Y_{\bar K}=V(q_\lambda,q_2,\ldots,q_m)
		\subset\PP^n_{\bar K}.
	$$
	Since $q_\lambda$ vanishes on $\operatorname{rad}(q_\lambda)$ and
	$$
		\dim\PP(\operatorname{rad}(q_\lambda))
		=\operatorname{corank}(q_\lambda)-1\ge m-1,
	$$
	the projective dimension theorem shows that
	$$
		\PP(\operatorname{rad}(q_\lambda))\cap V(q_2,\ldots,q_m)
	$$
	is nonempty: the $m-1$ restricted homogeneous equations cut out a subscheme of dimension at least
	$$
		\dim\PP(\operatorname{rad}(q_\lambda))-(m-1)\ge0.
	$$
	Choose $p=[v]$ in this intersection. Then $p\in Y_{\bar K}$ and $dq_\lambda|_v=0$. Hence the Jacobian matrix of $q_\lambda,q_2,\ldots,q_m$ at $v$ has rank at most $m-1$, contrary to the smoothness of $Y_{\bar K}$. Therefore
	$$
		\operatorname{corank}(q_\lambda)\le m-1.
	$$
	Smoothness and the complete intersection property are preserved by field extension, so the same arguments apply over every extension of $K$.
\end{proof}

\begin{lemma}\label{higher:quadric-family}
	Define the closed subscheme
	\begin{equation}\label{higher:incidence}
		Z= \left\{([\lambda],[w])\in\PP_K^{m-1}\times\PP(W) \ \middle|\
		      M(w)^{\mathsf t}\lambda=0,\quad \lambda^{\mathsf t}Q(w)=0
		      \right\}.
	\end{equation}
	Then $Z$ is integral, and the first projection $\pi:Z\longrightarrow\PP_K^{m-1}$ is flat and projective, with fibers that are quadric hypersurfaces of dimension $n-2m+1$. Its generic fiber is the smooth quadric obtained by hyperbolic reduction of the generic member of the linear system along $U$.
\end{lemma}

\begin{proof}
	Let $\lambda_1,\ldots,\lambda_m$ be homogeneous coordinates on $\PP_K^{m-1}$. Define the morphism of vector bundles
	$$
		\rho:W\otimes\cO_{\PP^{m-1}}
		\longrightarrow U^\vee\otimes\cO_{\PP^{m-1}}(1),
		\qquad \rho(w)(u)=\lambda^{\mathsf t}M(w)u.
	$$
	By Lemma~\ref{higher:rank}, $\rho$ is surjective on every fiber. With $\mathcal K=\ker\rho$, we obtain the exact sequence
	\begin{equation}\label{higher:kernel}
		0\longrightarrow\mathcal K\longrightarrow
		W\otimes\cO_{\PP^{m-1}}\xrightarrow{\rho}
		U^\vee\otimes\cO_{\PP^{m-1}}(1)\longrightarrow0.
	\end{equation}
	The kernel $\mathcal{K}$ has rank
	$$
		\operatorname{rank}\mathcal K = \dim W-\dim U =
		(n-m+2)-(m-1) = n-2m+3.
	$$
    The equations $M(w)^{\mathsf t}\lambda=0$ define the projective bundle
	$$
		P:=\PP_{\PP_K^{m-1}}(\mathcal K)
		\subset\PP_K^{m-1}\times\PP(W),
	$$
	where projectivization parametrizes one-dimensional subspaces. The bihomogeneous polynomial $\lambda^{\mathsf t}Q(w)$ defines a section
	$$
		s\in H^0\!\left(P,
		\left(\cO_{\PP_K^{m-1}}(1)\boxtimes\cO_{\PP(W)}(2)\right)|_P\right),
		\qquad Z=V_P(s).
	$$
	Fix $[\lambda]\in\PP_K^{m-1}$, let $\kappa$ be its residue field, and choose a nonzero representative $\lambda\in\kappa^m$. In the following argument, all vector spaces are extended to $\kappa$. The fiber of $\mathcal K$ at $[\lambda]$ is
	$$
		\mathcal K_{[\lambda]} = \{w\in W\mid M(w)^{\mathsf t}\lambda=0\}
		= W\cap U^{\perp_{q_\lambda}}.
	$$
	Since $q_\lambda|_U=0$, we have $U^{\perp_{q_\lambda}}=U\oplus\mathcal K_{[\lambda]}$. The quotient map therefore induces an isomorphism
	$$
		\mathcal K_{[\lambda]} \longrightarrow U^{\perp_{q_\lambda}}/U, \qquad  w\longmapsto[w].
	$$
	Under the quotient isomorphism $\mathcal K_{[\lambda]}\simeq U^{\perp_{q_\lambda}}/U$, the quadratic form
	$$
		\mathcal K_{[\lambda]}\longrightarrow\kappa,    \qquad  w\longmapsto\lambda^{\mathsf t}Q(w)
	$$
	is identified with the induced quadratic form
	$$
		\overline q_\lambda:U^{\perp_{q_\lambda}}/U\longrightarrow\kappa, \qquad [v]\longmapsto q_\lambda(v).
	$$
	By Lemma~\ref{higher:rank}, the map $w\mapsto b_{q_\lambda}(-,w)|_U$ from $W$ to $U^\vee$ is surjective. Choose a complement $A$ to $\mathcal K_{[\lambda]}$ in $W$. Then $\dim A=m-1$, and $b_{q_\lambda}$ induces a perfect pairing $U\times A\longrightarrow\kappa$. The restriction of $q_\lambda$ to $U\oplus A$ is nondegenerate. Indeed, let $v=u+a$, with $u\in U$ and $a\in A$, lie in $\operatorname{rad}(q_\lambda|_{U\oplus A})$. Since $U$ is totally isotropic, for every $u'\in U$ we have
	$$
		0=b_{q_\lambda}(v,u')
		=b_{q_\lambda}(u,u')+b_{q_\lambda}(a,u')
		=b_{q_\lambda}(a,u').
	$$
	Since $b_{q_\lambda}|_{U\times A}$ is a perfect pairing, $a=0$. Thus $v=u$, and for every $a'\in A$ we have
	$$
		0=b_{q_\lambda}(v,a')=b_{q_\lambda}(u,a').
	$$
	The nondegeneracy of $b_{q_\lambda}|_{U\times A}$ gives $u=0$. Thus $q_\lambda|_{U\oplus A}$ is nondegenerate. To describe $q_\lambda|_{U\oplus A}$ explicitly, choose bases $e_1,\ldots,e_{m-1}$ of $U$ and $a_1,\ldots,a_{m-1}$ of $A$ satisfying
	$$
		b_{q_\lambda}(e_i,a_j)=\delta_{ij}
		\qquad (1\leq i,j\leq m-1).
	$$
	For each $i\in\{1,\ldots,m-1\}$, put
	$$
		a_i' = a_i-\frac12 \sum_{j=1}^{m-1}b_{q_\lambda}(a_i,a_j)e_j.
	$$
	Then
	$$
		b_{q_\lambda}(a_i',a_j')=0, \qquad
		b_{q_\lambda}(e_i,a_j')=\delta_{ij}, \qquad
		q_\lambda(e_i)=q_\lambda(a_i')=0.
	$$
	Thus each $\kappa e_i+\kappa a_i'$ is a hyperbolic plane for $q_\lambda$, and these $m-1$ planes are pairwise orthogonal. Since $U\oplus A$ is nondegenerate, it is an orthogonal summand of $\kappa^{n+1}$. Moreover,
	$$
		U^{\perp_{q_\lambda}} = U\oplus(U\oplus A)^{\perp_{q_\lambda}},
	$$
	so the quotient map induces an isometry
	$$
		(U\oplus A)^{\perp_{q_\lambda}} \longrightarrow U^{\perp_{q_\lambda}}/U,
		\qquad v\longmapsto[v].
	$$
	With the hyperbolic plane $\mathbb H$ over $\kappa$ as in Section~\ref{pre:quadratic-forms}, we obtain the orthogonal decomposition
	\begin{equation}\label{higher:hyperbolic-decomposition}
		q_\lambda \simeq
		\underbrace{
		      \mathbb H\perp\cdots\perp\mathbb H
		}_{m-1\text{ copies}} \perp\overline q_\lambda.
	\end{equation}
	Each hyperbolic plane has dimension and rank two. Therefore,
	$$
		\begin{aligned}
			&\dim_\kappa\bigl(U^{\perp_{q_\lambda}}/U\bigr) =(n+1)-2(m-1) =n-2m+3,\\
			&\operatorname{rank}(\overline q_\lambda) =\operatorname{rank}(q_\lambda)-2(m-1),\\
			&\operatorname{corank}(\overline q_\lambda) =\operatorname{corank}(q_\lambda) \leq m-1.
		\end{aligned}
	$$
	It follows that
	$$
		\operatorname{rank}(\overline q_\lambda) \geq (n-2m+3)-(m-1)
		= n-3m+4 \geq2.
	$$
	Thus the restriction of $s$ to every geometric fiber $P_{[\lambda]}=\PP(\mathcal K_{[\lambda]})$ is nonzero. Each fiber $P_{[\lambda]}$ is integral, so $s|_{P_{[\lambda]}}$ is a nonzerodivisor. The zero scheme $Z_{[\lambda]}=V(s|_{P_{[\lambda]}})$ is a quadric hypersurface of dimension
	$$
		(\operatorname{rank}\mathcal K-1)-1 = n-2m+1.
	$$
	By \cite[Tag~062Y]{Stacks}, $Z$ is therefore a relative effective Cartier divisor in the sense of \cite[Tag~056P]{Stacks}. Hence $\pi:Z\to\PP_K^{m-1}$ is flat. The morphism $\pi$ is projective because it is the composite of the closed immersion
	$$
		Z\hookrightarrow\PP_K^{m-1}\times\PP(W)
	$$
	with the first projection.
	Let $\eta$ be the generic point of $\PP_K^{m-1}$. The form $q_\eta$ is nondegenerate by hypothesis, so \eqref{higher:hyperbolic-decomposition} shows that $\overline q_\eta$ is nondegenerate. Thus $Z_\eta=V(\overline q_\eta)$ is a smooth quadric. Since its dimension satisfies
	$$
		n-2m+1\geq m-1\geq1,
	$$
	$Z_\eta$ is geometrically integral. Moreover, $Z$ is an effective Cartier divisor in the smooth $K$-variety $P$, so $Z$ is Cohen--Macaulay. Flatness of $\pi$ implies that every irreducible component of $Z$ dominates $\PP_K^{m-1}$. Since $Z_\eta$ is integral, $Z$ is irreducible and generically reduced. The Cohen--Macaulay scheme $Z$ has no embedded points. Hence $Z$ is reduced and therefore integral.
\end{proof}

\begin{lemma}\label{higher:birational}
	Let
	$$
		\Omega= \{[w]\in\PP(W)\mid\operatorname{rank}M(w)=m-1\}, \qquad F(w)=\det[\,M(w)\mid Q(w)\,].
	$$
	Then $F$ is nonzero and $V(F)\cap\Omega$ is nonempty. The morphisms
	$$
		Y\setminus\PP(U)\longrightarrow\PP(W), \qquad [u:w]\longmapsto[w],
	$$
	and
	$$
		Z\longrightarrow\PP(W), \qquad ([\lambda],[w])\longmapsto[w],
	$$
	where $Z$ is the scheme of Lemma~\ref{higher:quadric-family}, induce isomorphisms from their respective inverse images of $V(F)\cap\Omega$ onto $V(F)\cap\Omega$. In particular, $Y$ and $Z$ are $K$-birational.
\end{lemma}

\begin{proof}
	\noindent\textbf{Step 1.} We prove that $\PP(W)\setminus\Omega$ has codimension at least $2$ in $\PP(W)$. It is enough to prove that
	$$
		\dim\bigl(\PP(W)\setminus\Omega\bigr)\leq n-m-1.
	$$
	Consider the incidence scheme
	$$
		S= \left\{ ([u],[w])\in\PP(U)\times\PP(W)\ \middle|\ M(w)u=0 \right\}.
	$$
	For each $[u]\in\PP(U)$, the fiber of the first projection $S\to\PP(U)$ is the projectivization of the kernel of
	$$
		W\longrightarrow K^m, \qquad w\longmapsto M(w)u,
	$$
	with scalar extension to the residue field of $[u]$ understood. By Lemma~\ref{higher:rank}, $w\mapsto M(w)u$ is surjective, so every fiber of $S\to\PP(U)$ has dimension
	$$
		\dim W-m-1=n-2m+1,
	$$
	and consequently
	$$
		\dim S=(m-2)+(n-2m+1)=n-m-1.
	$$
	The image of the second projection $S\to\PP(W)$ is $\PP(W)\setminus\Omega$. Thus $\dim(\PP(W)\setminus\Omega)\leq\dim S=n-m-1$. Since $\dim\PP(W)=n-m+1$, the complement of $\Omega$ has codimension at least $2$. In particular, $\Omega$ is nonempty.

	\medskip \noindent\textbf{Step 2.} We prove that $F\ne0$, that $V(F)\cap\Omega$ is nonempty, and that the projection
	$$
		\phi : Y\setminus\PP(U)\longrightarrow\PP(W),
		\qquad [u:w]\longmapsto[w],
	$$
	induces an isomorphism from the inverse image $\phi^{-1}(V(F)\cap\Omega)$ onto $V(F)\cap\Omega$. For $[w]\in\Omega$, the matrix $M(w)$ has rank $m-1$. Thus, for a fixed representative $w$, the equations
	$$
		M(w)u=-Q(w)
	$$
	have at most one solution. A solution exists if and only if $F(w)=\det[\,M(w)\mid Q(w)\,]=0$. Choose a set $I\subset\{1,\ldots,m\}$ of $m-1$ row indices, and write $M_I(w)$ and $Q_I(w)$ for the corresponding square submatrix of $M(w)$ and subvector of $Q(w)$. Set
	$$
		\Omega_I=\{\det M_I\ne0\}\subset\PP(W).
	$$
	On $\Omega_I$, the equations indexed by $I$ give
	$$
		u=-M_I(w)^{-1}Q_I(w).
	$$
	Let $j$ be the unique element of $\{1,\ldots,m\}\setminus I$. Substituting $u=-M_I(w)^{-1}Q_I(w)$ into the $j$-th equation
	$$
		\sum_{a=1}^{m-1}L_{ja}(w)u_a+Q_j(w)=0
	$$
	gives
	$$
		Q_j(w)-\sum_{a=1}^{m-1}L_{ja}(w)\bigl(M_I(w)^{-1}Q_I(w)\bigr)_a
		=(-1)^{j+m}\frac{F(w)}{\det M_I(w)}=0.
	$$
	Since $\det M_I(w)$ is nonzero on $\Omega_I$, this equation is equivalent to $F(w)=0$. Clearing the denominator with the adjugate matrix gives the morphism
	$$
		\begin{aligned}
			\psi_I :  V(F)\cap\Omega_I &\longrightarrow Y\setminus\PP(U),\\
			[w]&\longmapsto \bigl[-\operatorname{adj}(M_I(w))Q_I(w): \det(M_I(w))w \bigr].
		\end{aligned}
	$$
	The coordinates defining $\psi_I$ are homogeneous of degree $m$. They have no common zero on $\Omega_I$, since $\det M_I(w)w\ne0$. The morphism $\psi_I$ is inverse to $\phi$ over $V(F)\cap\Omega_I$. Suppose that $F=0$ as a polynomial. By Step 1, some $\Omega_I$ is nonempty. The morphism $\psi_I$ would then be defined on all of $\Omega_I$, so the projection $\phi$ would be dominant. This is impossible because
	$$
		\dim Y=n-m<n-m+1=\dim\PP(W).
	$$
	Hence $F\ne0$. The nonzero polynomial $F$ is homogeneous of degree $m+1$, so $V(F)$ is a hypersurface in $\PP(W)$. The complement $\PP(W)\setminus\Omega$ has codimension at least two, so
	$$
		V(F)\cap\Omega\ne\varnothing.
	$$
	The open subsets $\Omega_I$, for $|I|=m-1$, cover $\Omega$. Let $I,J\subset\{1,\ldots,m\}$ have cardinality $m-1$. For
	$$
		[w]\in V(F)\cap\Omega_I\cap\Omega_J,
	$$
	fix a representative $w$. Both $-M_I(w)^{-1}Q_I(w)$ and $-M_J(w)^{-1}Q_J(w)$ solve $M(w)u=-Q(w)$. Since $\operatorname{rank}M(w)=m-1$, the equation $M(w)u=-Q(w)$ has at most one solution. Therefore,
	$$
		\begin{aligned}
			\psi_I([w])
			&=[-M_I(w)^{-1}Q_I(w):w]\\
			&=[-M_J(w)^{-1}Q_J(w):w]
			=\psi_J([w]).
		\end{aligned}
	$$
	Thus $\psi_I$ and $\psi_J$ agree on the overlap $V(F)\cap\Omega_I\cap\Omega_J$, so the morphisms $\psi_I$ glue to an inverse of $\phi$ over $V(F)\cap\Omega$.

	\medskip \noindent\textbf{Step 3.} We prove that the projection
	$$
		\phi_Z:Z\longrightarrow\PP(W),
		\qquad ([\lambda],[w])\longmapsto[w],
	$$
	induces an isomorphism $\phi_Z^{-1}(V(F)\cap\Omega) \simeq V(F)\cap\Omega$. Together with Step 2, this shows that $Y$ and $Z$ are $K$-birational. For each $i\in\{1,\ldots,m\}$, let $M_{\widehat i}(w)$ be the square matrix obtained from $M(w)$ by deleting its $i$-th row, and put
	$$
		\lambda_i(w)=(-1)^{i+m}\det M_{\widehat i}(w), \qquad
		\lambda(w)=
		\begin{pmatrix}
			\lambda_1(w)\\
			\vdots\\
			\lambda_m(w)
		\end{pmatrix}.
	$$
	For $[w]\in\Omega$, the vector $\lambda(w)$ is nonzero and spans $\ker M(w)^{\mathsf t}$. Since each $\lambda_i(w)$ is homogeneous of degree $m-1$, the signed minors define a morphism
	$$
		\Omega\longrightarrow\PP_K^{m-1}, \qquad [w]\longmapsto[\lambda(w)].
	$$
	Expanding $F(w)=\det[\,M(w)\mid Q(w)\,]$ along the last column gives
	$$
		\lambda(w)^{\mathsf t}Q(w)=F(w).
	$$
	The identities $M(w)^{\mathsf t}\lambda(w)=0$ and $\lambda(w)^{\mathsf t}Q(w)=F(w)$ define a morphism
	$$
		\begin{aligned}
			\theta : V(F)\cap\Omega&\longrightarrow Z,\\
			[w]&\longmapsto([\lambda(w)],[w]).
		\end{aligned}
	$$
	On $\Omega_I$, let $j$ be the row index omitted from $I$. Then $\lambda_j(w)=(-1)^{j+m}\det M_I(w)$ is nonzero on $\Omega_I$. In the chart $\lambda_j=1$, the equations $M(w)^{\mathsf t}\lambda=0$ give
	$$
		\lambda=\frac{\lambda(w)}{\lambda_j(w)},
		\qquad \lambda^{\mathsf t}Q(w)=\frac{F(w)}{\lambda_j(w)}.
	$$
	Thus the equations of $Z$ over $\Omega_I$ reduce to $F(w)=0$. Hence $\theta$ is inverse to the restriction of $\phi_Z$ over $V(F)\cap\Omega$. The equations $M(w)u+Q(w)=0$ of $Y$ and \eqref{higher:incidence} of $Z$ give
	$$
		\begin{aligned}
			\phi^{-1}(\Omega)&=\phi^{-1}(V(F)\cap\Omega),\\
			\phi_Z^{-1}(\Omega)&=\phi_Z^{-1}(V(F)\cap\Omega).
		\end{aligned}
	$$
	By Steps 2 and 3, $\phi^{-1}(\Omega)$ and $\phi_Z^{-1}(\Omega)$ are isomorphic to the nonempty scheme $V(F)\cap\Omega$. Thus $\phi^{-1}(\Omega)$ and $\phi_Z^{-1}(\Omega)$ are nonempty open subsets of $Y$ and $Z$, respectively. Since $Y$ is integral by hypothesis and $Z$ is integral by Lemma~\ref{higher:quadric-family}, these open subsets are dense. The restrictions of $\phi$ and $\phi_Z$ to these dense open subsets give inverse birational maps
	$$
		Y\dashrightarrow Z, \qquad Z\dashrightarrow Y.
	$$
\end{proof}

\section{An odd-degree class on the conic fibration}\label{app:nonvacuous}
We prove the two facts used in Proposition~\ref{pes:nonvacuous}. By \cite[Proposition~5.2]{BS21}, choose a smooth nonempty Zariski open subset $U'\subset U$ that the period map identifies with an open subset of the normalization $\widetilde{\mathcal D}_{24}$, over the locus where $\widetilde{\mathcal D}_{24}\to\mathcal D_{24}$ is an isomorphism.

The incidence variety of pairs $([\sigma],V_1\subset V_6)$ satisfying $\sigma(V_1,V_6,V_{10})=0$ is a projective bundle over $\operatorname{Fl}(1,6,V_{10})$; see \cite[Remark~5.3]{BS21}. Hence it is smooth and irreducible. On a nonempty open subset, the constructions of \cite[\S5.1]{BS21} give smooth projective families $X_6^\sigma$, $D$, and $S$, with the inclusion $i$ and the conic fibration $\pi$. The natural map from this open subset to $\mathcal D^{1,6,10}$ is dominant. By generic smoothness, after restricting this open subset and shrinking $U'$, we obtain a smooth surjective map $f:T\to U'$ with $T$ smooth and connected. We use singular cohomology on the associated analytic spaces.

\begin{lemma}\label{pes:odd-pairing}
	For every $\sigma\in U'$, let $i:D\hookrightarrow X_6^\sigma$ be the inclusion and let $\ell\in H_2(D,\ZZ)$ be the class of a fiber of $\pi$. Then
	\begin{equation}\label{pes:fiber-degree}
		i^*C\cdot\ell=q(C,[D])
	\end{equation}
	for every $C\in H^2(X_6^\sigma,\ZZ)$.
\end{lemma}

\begin{proof}
	Take a very general point $\sigma_0\in U'$. Then $\operatorname{NS}(X_6^{\sigma_0})\otimes\QQ=\QQ h+\QQ[D]$, where $h$ is the polarization class; see \cite[\S5.2]{BS21}. As in the proof of \cite[Lemma~5.11]{BS21}, the class $y\in H^2(X_6^{\sigma_0},\QQ)$ defined by $q(C,y)=i^*C\cdot\ell$ for all $C\in H^2(X_6^{\sigma_0},\ZZ)$ has type $(1,1)$, because $i_*\ell$ is algebraic. Thus $y=ah+b[D]$ for some $a,b\in\QQ$. The fiber degrees and the matrix in \cite[Lemma~5.9]{BS21} give
	$$
		q(h,y)=22a+2b=2,\qquad q([D],y)=2a-2b=-2.
	$$
	Thus $y=[D]$, proving~\eqref{pes:fiber-degree} at $\sigma_0$. We choose $t_0\in T$ over $\sigma_0$. For $t\in T$, write $i_t:D_t\hookrightarrow X_6^{f(t)}$ for the divisor inclusion and $\ell_t$ for the conic fiber class. A path from $t_0$ to $t$ gives parallel transport
	$$
		\tau:H^2(X_6^{\sigma_0},\ZZ)\xrightarrow{\sim}H^2(X_6^{f(t)},\ZZ)
	$$
	in the cohomology local system; see \cite[\S\S9.1.1 and 9.2.1]{Voi02} and \cite[Definition~1.1]{Mar11}. Transport preserves $q$ and $[D]$. It also preserves $i_*\ell$: under Poincar\'e duality, this class is $i_*\pi^*[s]$, where $[s]\in H^4(S,\ZZ)$ is the point class, and pullback and Gysin maps commute with transport. Hence, for every $C\in H^2(X_6^{\sigma_0},\ZZ)$,
	$$
		\begin{aligned}
			i_t^*\tau(C)\cdot\ell_t
			&=i_{t_0}^*C\cdot\ell_{t_0}\\
			&=q(C,[D_{t_0}])=q(\tau(C),[D_t]).
		\end{aligned}
	$$
	Since $\tau$ is an isomorphism and $f$ is surjective, this proves~\eqref{pes:fiber-degree} for every $\sigma\in U'$.
\end{proof}

\begin{lemma}\label{pes:odd-class}
	There exist $\sigma\in U'$ and $C\in\operatorname{NS}(X_6^\sigma)$ such that $q(C,[D])=1$.
\end{lemma}

\begin{proof}
	Choose a lattice isometry $H^2(X_6^{\sigma_0},\ZZ)\simeq\Lambda$ at a very general point $\sigma_0\in U'$, and use $h,[D]$ for the images of the two divisor classes. By \cite[Lemma~5.9 and Corollary~5.10]{BS21}, the lattice $K=\ZZ h+\ZZ[D]$ is saturated, its Gram matrix is
	$$
		\begin{pmatrix}22&2\\2&-2\end{pmatrix},
	$$
	and there is $C_0\in\Lambda$ with $q(C_0,[D])=1$. Put $W=K^\perp\otimes\mathbb R$, which has signature $(2,19)$. The period domain in $W$ parametrizes oriented positive two-planes. Let $\Omega$ be the connected component presenting $\widetilde{\mathcal D}_{24}$, and let $O\subset\Omega$ be the inverse image of $U'$. The arithmetic group defining this quotient preserves $h$ and $K$. It also fixes $[D]$, since $[D]$ is the unique vector in $K$ of square $-2$ and pairing $2$ with $h$. With the identification induced by the chosen isometry, $[D]$ represents the geometric divisor throughout $O$. We choose $P\in O$ and a nonzero vector $v\in P^\perp\cap W$. Then $q(v)<0$. Since $K^\perp$ is a full lattice in $W$, rounding coordinates in an integral basis gives vectors $v_n\in K^\perp$ with $v_n/n\to v$. Set $C_n=C_0+v_n$, and let $w_n$ be the orthogonal projection of $C_n/n$ onto $W$. Then
	$$
		q(C_n,[D])=1,\qquad w_n\longrightarrow v,\qquad q(w_n)<0\quad(n\gg0).
	$$
	For sufficiently large $n$, define the two-plane
	$$
		P_n=\left\{x-\frac{q(x,w_n)}{q(w_n)}w_n:x\in P\right\},
	$$
	with the orientation induced by $P$. Since $q(x,v)=0$ for $x\in P$, the displayed map tends to the identity on $P$. Hence $P_n$ is positive and $P_n\to P$. As $O$ is open, $P_n\in O$ for all sufficiently large $n$. By construction, $P_n\subset W\cap C_n^\perp$. At the corresponding point $\sigma_n\in U'$, the integral class represented by $C_n$ therefore has type $(1,1)$ and pairs to one with $[D]$. The Lefschetz $(1,1)$ theorem places this class in $\operatorname{NS}(X_6^{\sigma_n})$.
\end{proof}

\bibliographystyle{amsalpha}
\bibliography{reference}

\end{document}